\documentclass[final,onefignum,onetabnum]{siamonline250106}

\usepackage{amssymb}
\usepackage{hyperref}
\usepackage{url}
\usepackage{graphicx}
\usepackage{chngcntr}
\usepackage{pgfplots}
\usepackage{wrapfig}
\usepackage{bm}
\usepackage{booktabs}
\usepackage{multirow}
\usepackage{multicol}
\usepackage{array} 

\usepackage{algorithm}        
\usepackage{algorithmicx}
\usepackage[noend]{algpseudocode} 
\usepgfplotslibrary{fillbetween}
\usetikzlibrary{arrows.meta}
\usetikzlibrary {shapes.geometric}
\tikzset{>=latex}
\definecolor{blue}{rgb}{0,0.45,0.74}
\definecolor{purple}{rgb}{0.76,0.48,0.82}
\definecolor{yellow}{rgb}{0.93,0.69,0.13}
\definecolor{green}{rgb}{0.47,0.67,0.19}
\definecolor{red}{rgb}{0.8500,0.33,0.1}
\usepackage{xcolor}

\pgfmathdeclarefunction{relu}{1}{%
    \pgfmathparse{max(#1,0)}%
}
\ifpdf
  \DeclareGraphicsExtensions{.eps,.pdf,.png,.jpg}
\else
  \DeclareGraphicsExtensions{.eps}
\fi

\usepackage{enumitem}
\setlist[enumerate]{leftmargin=.5in}
\setlist[itemize]{leftmargin=.5in}

\newsiamremark{remark}{Remark}
\newsiamremark{hypothesis}{Hypothesis}
\crefname{hypothesis}{Hypothesis}{Hypotheses}
\newsiamthm{claim}{Claim}

\headers{Primal-dual image reconstruction with ICNNs}{M. J. Ehrhardt, S. Mukherjee, and H. S. Wong}

\title{A primal-dual algorithm for image reconstruction with input-convex neural network regularizers\thanks{Submitted to the editors April 15, 2025.
\funding{Matthias J. Ehrhardt acknowledges support from the EPSRC (EP/S026045/1, EP/T026693/1, EP/V026259/1). Hok Shing Wong acknowledges support from the project EP/V026259/1. }}}

\author{Matthias J. Ehrhardt\thanks{Department of Mathematical Sciences, University of Bath, UK 
(\email{m.ehrhardt@bath.ac.uk},\email{hsw43@bath.ac.uk}).}
\and Subhadip Mukherjee\thanks{Department of Electronics \& Electrical Communication Engineering, Indian Institute of Technology (IIT), Kharagpur, India 
  (\email{smukherjee@ece.iitkgp.ac.in}).}
\and Hok Shing Wong\footnotemark[2]}

\usepackage{amsopn}

\ifpdf
\hypersetup{
  pdftitle={Primal-dual image reconstruction with ICNNs},
  pdfauthor={M. J. Ehrhardt, S. Mukherjee, and H. S. Wong}
}
\fi

\newcommand{\rv}[1]{\textcolor{black}{#1}}

\begin{document}

\maketitle

\begin{abstract}
We address the optimization problem in a data-driven variational reconstruction framework, where the regularizer is parameterized by an input-convex neural network (ICNN). While gradient-based methods are commonly used to solve such problems, they struggle to effectively handle non-smooth problems which often leads to slow convergence. Moreover, the nested structure of the neural network complicates the application of standard non-smooth optimization techniques, such as proximal algorithms. To overcome these challenges, we reformulate the problem and eliminate the network's nested structure. By relating this reformulation to epigraphical projections of the activation functions, we transform the problem into a convex optimization problem that can be efficiently solved using a primal-dual algorithm. We also prove that this reformulation is equivalent to the original variational problem. Through experiments on several imaging tasks, we show that the proposed approach not only outperforms subgradient methods and even accelerated methods in the smooth setting, but also facilitates the training of the regularizer itself.
\end{abstract}



\section{Introduction}
\label{sec:intro}
Image restoration focuses on reconstructing high-quality images from degraded, low-quality versions that often result from issues during image acquisition and transmission. This includes tasks such as image denoising, image deblurring, image inpainting and computer tomography (CT) reconstruction. The measurement process is typically modeled as $\mathbf{y} = \mathbf{Ax}+\boldsymbol{\epsilon}$, where $\mathbf{A}$ simulates the physics in the measurement process and $\boldsymbol{\epsilon}$ denotes the measurement noise. One then seeks to recover the unknown image $\mathbf{x}$ from the noisy measurement $\mathbf{y}$. To mitigate the ill-possedness of the inverse problem, the classical variational reconstruction framework incorporates prior information about plausible reconstructions through a regularizer:
\begin{equation}\label{var-prob}\tag{P}
    \min_{\mathbf{x}} D(\mathbf{Ax},\mathbf{y})+\gamma R_{\boldsymbol{\theta}}(\mathbf{x}),
\end{equation}
where $D$ is the data fidelity. The regularizer $R_{\boldsymbol{\theta}}$ can be parametrized, with $\mathbf{\theta}$ denoting its parameters. The trade-off between data fidelity and regularizer is controlled by the positive regularization parameter $\gamma$. The reconstruction is obtained by solving the minimization problem (\ref{var-prob}).

\subsection{Overview of learned regularizers} 
Traditional methods often utilize hand-crafted regularizers, such as total variation (TV) \cite{rudin1992nonlinear}, total generalized variation (TGV) \cite{bredies2010total} and sparsity promoting regularizer \cite{daubechies2004iterative}. \rv{Under mild assumptions, mathematical properties including convergence and stability can be established for these regularizers. See \cite{benning2018modern,scherzer2009variational} for detailed analysis.} In recent years, data-driven approaches for inverse problems have gained increasing interest. For instance, \cite{chen2017low,jin2017deep,kang2017deep} propose learning end-to-end neural networks that post-process analytical reconstructions. Another prominent strategy involves unrolling methods \cite{adler2018learned,kobler2017variational,meinhardt2017learning,yang2016deep}, which integrate neural network modules into iterative optimization algorithms based on the variational framework. 

\rv{A complementary approach seeks to combine data-driven methods with classical priors by learning explicit regularizers \cite{aharon2006k,chen2014insights,kunisch2013bilevel,xu2012low}, with the aim of obtaining the best of both worlds: the mathematical foundation and interpretability of handcrafted regularizers together with the reconstruction quality of data-driven models. For example, \cite{goujon2023neural, goujon2024learning} learns the convolutional kernels and potentials in Field of Experts (FoE) models \cite{roth2009fields}. This is also extended to parametrizing them with neural networks \cite{kobler2020total,li2020nett,lunz2018adversarial,mukherjee2020learned}, and embedding them within variational reconstruction frameworks. Similar ideas have also been applied to patch-based methods, where the regularizer is learned from local information \cite{altekruger2023patchnr,prost2021learning,zoran2011learning}. For a more comprehensive review of deep learning techniques for designing regularizers, we refer the reader to \cite{arridge2019solving,ongie2020deep} and references therein.}

\subsection{Existing methods for solving variational problems}
Numerous efforts have been made in the literature to study algorithms for optimizing convex functions, in particular for variational reconstruction. Gradient methods are often applied to general smooth convex problems \cite{boyd2004convex} and can be extended to subgradient methods for non-smooth problems \cite{boyd2003subgradient}. Another essential component for non-smooth problems is the proximal operator \cite{parikh2014proximal}. In particular, primal-dual methods have been extensively studied for non-smooth handcrafted regularizers such as TV \cite{chambolle2011first,chambolle2016ergodic,yan2018new,zhu2008efficient}. However, due to the nested structure of neural networks, computing the proximal operator for neural networks is often impractical. Therefore, to perform variational reconstruction with neural network-parameterized regularizers \rv{including ICNNs}, subgradient methods are commonly applied, where subgradients are computed via backpropagation \cite{mukherjee2020learned}. Despite the simplicity of this approach, challenges arise due to non-smoothness, step-size selection. On the other hand, the idea of removing nested structure of neural network has been explored \cite{askari2018lifted,carreira2014distributed,li2019lifted,taylor2016training,wang2023lifted,zhang2017convergent} in the context of training. These unconventional training approaches introduce auxiliary variables corresponding to layer-wise activations. The resulting formulations are typically relaxed by adding penalties to the induced equality constraints. However, the problem remains non-convex, and the minimizers are altered as a result of these relaxations.

\rv{For smooth problems, accelerated methods such as Polyak's heavy-ball method \cite{polyak1964some} and Nesterov's method \cite{nesterov1983method} have been proposed. For convex but non-smooth problem, FISTA \cite{beck2009fast} provides accelerated convergence. In \cite{goujon2024learning}, a safeguard accelerated gradient descent was introduced to solve the associated variational problem of a weakly convex regularizer. NMAPG \cite{li2015accelerated} extends to the non-convex and non-smooth setting and was applied for non-convex regularizers in \cite{hertrich2025learning}. However, for non-smooth regularizers, these proximal-based accelerations remain impractical when the regularizer does not admit an efficient proximal mapping.} 

\subsection{Contributions}
Primal-dual algorithms have been successfully applied to classical variational problems, providing fast reconstruction methods. Motivated by their flexibility and practicality, we aim to exploit both the inherent convex nature and the architecture of the neural network to devise optimization algorithm for solving the variational problem. Our contributions are as follows: 
\begin{itemize}
    \item We introduce a more general architecture than ICNN. To address the non-smoothness and nested structure, we propose a novel reformulation of the variational problem. We prove that this reformulation is both convex and equivalent to the original variational problem.
    \item We apply this novel convex reformulation to setting where the regularizer is parameterized by an ICNN, solving the associated variational problem using a primal-dual algorithm. Additionally, we design a step-size scheme tailored specifically to our formulation.
    \item We implement the proposed framework for image restoration tasks such as denoising, inpainting, and CT reconstruction. Our results demonstrate that the proposed method is superior to subgradient methods \rv{and even accelerated method for smooth regularizers}, achieving faster and more stable reconstruction.
    \item \rv{We further incorporate the proposed method into the training pipeline by employing it as the lower-level solver in a bilevel setting for learning the regularizer, thereby demonstrating its effectiveness not only for reconstruction but also within the training process.}
\end{itemize}

\section{Background}
In this paper, we focus on solving (\ref{var-prob}), where the regularizer $R_{\boldsymbol{\theta}}$ is parameterized by an ICNN. A major advantage of a convex setting over a non-convex one is the ability to compute a global optimum independent of initialization, allowing one to leverage the well-established theory of convex optimization with guaranteed convergence to efficiently solve (\ref{var-prob}). \rv{In what follows, we provide further details on ICNN-based regularizers.}

For $\mathbf{x,y}\in\mathbb{R}^n$, we denote $\mathbf{x}\leq \mathbf{y}$ if $\mathbf{x}_i\leq \mathbf{y}_i$ for $i=1,\dots,n$. To handle general activations,  we call a function $f:\mathbb{R}^n\rightarrow\mathbb{R}^m$ convex if $f(\alpha \mathbf{x}+(1-\alpha)\mathbf{y})\leq\alpha f(\mathbf{x})+(1-\alpha)f(\mathbf{y})$ for every $\mathbf{x,y}\in\mathbb{R}^n$ and $\alpha\in[0,1]$. $f$ is called non-decreasing if $f(\mathbf{x})\leq f(\mathbf{y})$ for $\mathbf{x}\leq \mathbf{y}$.
\subsection{Architecture of ICNNs}
\label{sec:icnn-recall}
In \cite{amos2017input}, a $L$-layered ICNN is defined by the following architecture:
\begin{equation}\label{icnn-arch}
\tag{EQ}
\begin{aligned}
    \mathbf{z}_{1} &= h_{1}(\mathbf{V}_{0}\mathbf{x}+\mathbf{b}_0),\\
    \mathbf{z}_{i+1} &= h_{i+1}(\mathbf{V}_{i}\mathbf{x}+\mathbf{W}_{i}\mathbf{z}_i+\mathbf{b}_i), \ i=1,\dots,L-2,\\
    R_{\theta}(\mathbf{x})&:=h_{L}(\mathbf{V}_{L-1}\mathbf{x}+\mathbf{W}_{L-1}\mathbf{z}_{L-1}+\mathbf{b}_{L-1}),
\end{aligned}
\end{equation}
where $\mathbf{V}_{i},\mathbf{W}_{i}$ are linear operators, which could represent various neural network components, such as fully connected layers, convolution layers and average pooling layers. Here $\boldsymbol{\theta}=\{\mathbf{V}_{i},\mathbf{W}_{i},\mathbf{b}_i\}$ represents the collection of all trainable parameters of the ICNN. The functions $h_i$ are non-linear activations. 

Utilizing the fact that both non-negative sums
of convex functions and composition of a convex and a convex non-decreasing function are convex, the convexity of $R_{\boldsymbol{\theta}}$ with respect to the input $\mathbf{x}$ can be guaranteed by imposing that the weights $\mathbf{W}_{i}$ are non-negative and $h_i$ are convex, non-decreasing. 

\subsection{Training Methods}
\rv{In this section, we review training methods for learned regularizers, including those for ICNNs. A variety of approaches have been proposed in the literature, including adversarial regularization \cite{lunz2018adversarial,mukherjee2020learned} and bilevel learning \cite{calatroni2017bilevel,crockett2022bilevel}. While other strategies, such as Bayesian methods \cite{altekruger2023patchnr,ardizzone2018analyzing} have also been explored, we focus on the first two methods, as these are the approaches employed in this work.}

\subsubsection{Adversarial regularization framework}
\rv{The idea of adversarial regularization (AR) was first introduced in \cite{lunz2018adversarial}. It was applied to learn parameters of ICNN-based regularizers in \cite{mukherjee2020learned}. The goal of this scheme is to train the regularizer as a classifier, so that it outputs low values when provided with true images and higher values for task-dependent unregularized reconstructions. Given training images $\mathbf{x}_i,i=1,\dots,N_1$, and noisy measurements $\mathbf{y}_j,j=1,\dots,N_2$, the following training problem was considered:}
\rv{
\begin{equation}
\label{eq:AR_loss}
    \min_{\boldsymbol{\theta}} \frac{1}{N_1}
    \sum_{i=1}^{N_1} R_{\boldsymbol{\theta}}(\mathbf{x}_i)
    -
    \frac{1}{N_2}
    \sum_{j=1}^{N_2} R_{\boldsymbol{\theta}}(\mathbf{A}^{\dagger} \mathbf{y}_j) + \lambda_{GP} \mathbb{E}_{\mathbf{x}}\bigl[\bigl(\left\|\nabla R_{\boldsymbol{\theta}}(\mathbf{x})\right\|-1\bigr)_{+}^2\bigr],
\end{equation}
}
\rv{here $\mathbf{A}^{\dagger}$ denotes a (regularized) pseudo-inverse, and the expectation is taken over points along lines joining the samples $x_i$ and $\mathbf{A}^{\dagger}\mathbf{y}_j$. This term enforce approximate 1-Lipschitz $R_{\boldsymbol{\theta}}$, which is motivated by the connection to the 1-Wasserstein loss. Note that this setup can be applied in a semi-supervised setting, as the images and the noisy measurements are not necessarily paired.}

\rv{Motivated by the simplicity and relatively fast training of this method, we primarily adopt the adversarial framework to train the ICNN-based regularizers considered in this work.}

\subsubsection{Bilevel Learning}
\rv{While adversarial training is simple and often computationally fast, it does not explicitly account for reconstruction quality during training. Moreover, selecting an appropriate regularization parameter for adversarially trained regularizers is nontrivial. Bilevel learning (BL) provides a natural way to address these issues by aligning the training objective with reconstruction performance. Given training data pairs $(\mathbf{x}_i,\mathbf{y}_i)$, $i=1,...,N$, the regularizer can be trained by considering the following nested optimization problem:}
\rv{
\begin{equation}\label{eq:bilv}
    \begin{aligned}
        &\min_{\boldsymbol{\theta}} \left\{\mathcal{L}(\boldsymbol{\theta}) = \frac{1}{N}\sum_{i=1}^N
    \ell(\hat{\mathbf{x}}_{\mathbf{y}_i}(\boldsymbol{\theta}))\right\}\\
    &\hat{\mathbf{x}}_{\mathbf{y}_i}(\boldsymbol{\theta}) = \arg\min\limits_{\mathbf{x}}\left\{\mathcal{J}_{\mathbf{y}_i}(\mathbf{x}; \boldsymbol{\theta}) = D(\mathbf{Ax,y_i}) + R_{\boldsymbol{\theta}}(\mathbf{x}) \right\},
    \end{aligned}
\end{equation}
}

\rv{Here, the lower-level problem corresponds to the standard variational reconstruction, while the upper-level problem tunes the regularizer parameters $\boldsymbol{\theta}$ to minimize a reconstruction loss $\ell$, such as the $\ell_2$ error. Solving this nested optimization problem often requires computing the gradient of $\mathcal{L}$ with respect to $\boldsymbol{\theta}$, also known as the hypergradient. For a single data point, this can be computed using chain rule as $\nabla\mathcal{L}(\boldsymbol{\theta})=[\hat{\mathbf{x}}'_{\mathbf{y}}(\boldsymbol{\theta})]^T\nabla\ell(\hat{\mathbf{x}}_{\mathbf{y}}(\boldsymbol{\theta}))$. The central element for hypergradient computation is  $\hat{\mathbf{x}}'_{\mathbf{y}}(\boldsymbol{\theta})$. Various approaches for computating this have been proposed, including implicit differentiation (IFT) \cite{bengio2000gradient, grazzi2020iteration, samuel2009learning}, which relies on differentiating the optimality conditions of the lower-level problem but it involves Hessian inversion. Unrolling methods \cite{mehmood2020automatic, ochs2016techniques} on the other hand apply automatic differentiation to backpropagate through the iterations of the lower-level solver. However, this can be are memory-intensive as it requires storing the intermediate interates. Alternatively,  the Jacobian-free backpropagation (JFB) \cite{bolte2023one, fung2022jfb} reduces memory usage by only backpropagating through a truncated number of optimization steps. For more details on bilevel optimization, see \cite{ji2021bilevel,zucchet2022beyond}.}
\section{Constrained Convex Reformulation}
\label{sec:reform}
\rv{Rather than resorting to general (sub)gradient-based method, we propose a tailored approach for solving the variational problem associated with ICNN-based regularizers. The central idea is to resolve the nested structure of neural networks. By reformulating the problem, we replace the need for the intractable proximal mapping of the regularizer with simpler, computable operations. To introduce our proposed reformulation, we first consider a more general nested structure than the specific ICNN architecture for the functional $R_{\boldsymbol{\theta}}$:}
\begin{equation}\label{gen-eq}
\tag{EQ-G}
\begin{aligned}
    \mathbf{z}_{1} &= \phi_{1}(\mathbf{x}),\\
    \mathbf{z}_{i+1} &= \phi_{i+1}(\mathbf{x},\boldsymbol{\omega}_i)\text{ for }i=1,\dots,L-2, \text{ with }\boldsymbol{\omega}_i=(\mathbf{z}_1,\dots,\mathbf{z}_i),\\ 
    R_{\boldsymbol{\theta}}(\mathbf{x})&=\phi_L(\mathbf{x},\boldsymbol{\omega}_{L-1}).
\end{aligned}
\end{equation}

We make the following assumption on the activation functions.

\textbf{Assumption 1. } 
$\phi_i$ are convex for $i=1,\dots,L$, and $\phi_i^\mathbf{x}$ are non-decreasing for $i=2,\dots,L$, where $\phi_i^\mathbf{x}(\boldsymbol{\omega}_{i-1})=\phi_i(\mathbf{x},\boldsymbol{\omega}_{i-1})$.

\begin{proposition}\label{cvx}
Under Assumption 1, $R_{\boldsymbol{\theta}}$ defined by (\ref{gen-eq}) is convex with respect to $\mathbf{x}$.
\end{proposition}
\begin{proof}
Consider $\bar{\mathbf{x}},\tilde{\mathbf{x}}$, and $\lambda\in[0,1]$. Then, 
\begin{equation*}
    \mathbf{z}^{\lambda}_1:=\phi_1(\lambda\bar{\mathbf{x}}+(1-\lambda)\tilde{\mathbf{x}})\leq\lambda\phi_1(\bar{\mathbf{x}})+(1-\lambda)\phi_1(\tilde{\mathbf{x}})=:\lambda\bar{\mathbf{z}}_1+(1-\lambda)\tilde{\mathbf{z}}_1,
\end{equation*}
where the inequality is due the convexity of $\phi_1$. Since $\phi^\mathbf{x}_2$ is non-decreasing, we have:
\begin{equation}
\begin{aligned}
    \mathbf{z}^{\lambda}_2&:=\phi_2(\lambda\bar{\mathbf{x}}+(1-\lambda)\tilde{\mathbf{x}},\boldsymbol{\omega}^{\lambda}_1)\\
    &\leq\phi_2(\lambda\bar{\mathbf{x}}+(1-\lambda)\tilde{\mathbf{x}},\lambda\bar{\boldsymbol{\omega}}_1+(1-\lambda)\tilde{\boldsymbol{\omega}}_1)\\
    &\leq\lambda\bar{\mathbf{z}}_2+(1-\lambda)\tilde{\mathbf{z}}_2,
\end{aligned}
\end{equation}
where the second inequality follows from the convexity of $\phi_2$. Using similar argument, we have: 
\begin{equation}
    \boldsymbol{\omega}^{\lambda}_{i}\leq\lambda\bar{\boldsymbol{\omega}}_{i}+(1-\lambda)\tilde{\boldsymbol{\omega}}_{i}, \text{ for }i=2,\dots,L-1,
\end{equation}
where $\boldsymbol{\omega}_i$ are defined as $(\mathbf{z}_1,\dots,\mathbf{z}_i)$ and $\boldsymbol{\omega}_i^{\lambda}:=\phi_i(\lambda\bar{\mathbf{x}}+(1-\lambda)\tilde{\mathbf{x}},\boldsymbol{\omega}_{i-1}^{\lambda})$. 
In particular:
\begin{equation}
\begin{aligned}
    R_{\boldsymbol{\theta}}(\lambda\bar{\mathbf{x}}+(1-\lambda)\tilde{\mathbf{x}})&=\phi_L(\lambda\bar{\mathbf{x}}+(1-\lambda)\tilde{\mathbf{x}},\boldsymbol{\omega}^{\lambda}_{L-1})\\
    &\leq\phi_L(\lambda\bar{\mathbf{x}}+(1-\lambda)\tilde{\mathbf{x}},\lambda\bar{\boldsymbol{\omega}}_{L-1}+(1-\lambda)\tilde{\boldsymbol{\omega}}_{L-1})\\
    &\leq\lambda\phi_L(\bar{\mathbf{x}},\bar{\boldsymbol{\omega}}_{L-1})+(1-\lambda)\phi_L(\tilde{\mathbf{x}},\tilde{\boldsymbol{\omega}}_{L-1}),
\end{aligned}
\end{equation}
where the first inequality holds since $\phi^\mathbf{x}_L$ is non-decreasing and the second inequality is due to the convexity of $\phi_L$. Hence, $R$ is convex with respect to $\mathbf{x}$.
\end{proof}

Note that the ICNN architecture given by (\ref{icnn-arch}) is a special case of the above structure, with $\phi_{i+1}(\mathbf{x},\boldsymbol{\omega}_i)=h_{i+1}(\mathbf{V}_i\mathbf{x}+\mathbf{W}_i\mathbf{z}_i+\mathbf{b}_i)$. In particular, $\mathbf{W}_i$ being non-negative and $h_{i+1}$ being non-decreasing imply that $\phi_{i+1}^\mathbf{x}$ is non-decreasing. Hence, $R_{\boldsymbol{\theta}}$ parametrized as in (\ref{icnn-arch}) is indeed convex. We also relax the condition on $h_1$ to be merely convex, rather than both convex and non-decreasing, as in \cite{amos2017input}. With the above framework, we could also consider a residual architecture, where $\phi_{i+1}(\mathbf{x},\boldsymbol{\omega}_i)=\mathbf{z}_i+h_{i+1}(\mathbf{V}_i\mathbf{x}+\mathbf{W}_i\mathbf{z}_i+\mathbf{b}_i)$.

The main objective of this paper is to minimize a functional $R_{\boldsymbol{\theta}}$ with the above structure efficiently. The first step of the proposed approach involves removing the nested structure of the problem. Given $R_{\boldsymbol{\theta}}$ as defined in (\ref{gen-eq}), the problem (\ref{var-prob}) is equivalent \cite{carreira2014distributed} to: 
\begin{equation}\label{var-g}
\min_{\mathbf{x},\boldsymbol{\omega}_{L-1}}D(\mathbf{Ax,y})+\gamma\phi_L(\mathbf{x},\boldsymbol{\omega}_{L-1})\text{  subject to $\boldsymbol{\omega}_{L-1}$ satisfying (\ref{gen-eq})}. 
\end{equation}

However, the above reformulation is in general not convex as $\phi_i$ could be non-linear.

\textbf{Example.}
To illustrate the non-convexity of (\ref{var-g}), consider a simple 1D example. Here, we define $R_{\boldsymbol{\theta}}(x)=\exp(x+\max(x,0))$ and a data fidelity $D(x,y)=\frac{1}{2}(x-y)^2$. Then reformulation (\ref{var-g}) can be written as: 
\begin{equation*}
        \min_{x,z}\frac{1}{2}(x-y)^2+\exp(x+z)\text{ 
 subject to }z=\max(x,0).
\end{equation*}

\begin{wrapfigure}{R}{4.5cm}
\vspace*{-.5cm}
\begin{tikzpicture}
\begin{axis}[
  axis equal image,
  axis lines=middle,
  xlabel={$x$},
  ylabel={$f(x)$},
  xtick={\empty},
  ytick={\empty},
  domain=-1.25:1.25,
  ymin=-0.12,
  xmin=-1.5,
  xmax=1.25,
  clip=false,
  height=5cm
] 
\addplot [red,name path=B,samples=100] plot {relu(x)};
\addplot [no marks,draw=white,name path=C] coordinates 
  {(-1.25,{relu(-1.25)+0.5}) (1.25,{relu(1.25)+0.1})};
\makeatother
\addplot[draw=white,top color=gray!80!green!05,bottom color=gray!90!green!80] 
  fill between[of=B and C,
  soft clip={domain=-1.25:1.25}
  ];
\addplot [color=blue,only marks,mark=x,mark options={solid}]
    coordinates{(-1,0) (1,1) (0,0.5)};
\addplot [color=blue,dashed]
    coordinates{(-1,0) (1,1)};
\node at (axis cs:-1,-0.1) {$\mathbf{w}_1$};
\node at (axis cs:1,0.9) {$\mathbf{w}_2$};
\node at (axis cs:0,0.4) {$0.5\mathbf{w}_1+0.5\mathbf{w}_2$};
\end{axis}
\end{tikzpicture}
\vspace*{-.3cm}
    \caption{Example illustrating the non-convexity of (\ref{var-g}).}
    \label{fig:non-cvx}
\end{wrapfigure}
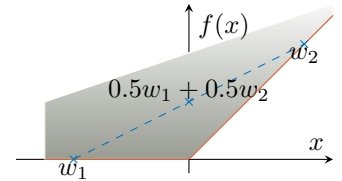

Here $\mathbf{w}_1=(-1,0),\mathbf{w}_2=(1,1)$ are both feasible but $0.5\mathbf{w}_1+0.5\mathbf{w}_2=(0,0.5)$ is not. Hence, the above problem is non-convex despite that the objective is convex. This is due to the fact that the graph (red curve) of the $\max$ function is not a convex set. However, $0.5\mathbf{w}_1+0.5\mathbf{w}_2$ belongs to the shaded region given by $\{(x,z)|z\geq \max(x,0)\}$, which is the epigraph of $\max$. In fact, epigraphs can represent a large family of non-linear constraints which are effective in inverse problems. Epigraphical projections were applied in  \cite{chierchia2015epigraphical} to solve classes of constrained convex optimization problems.

This motivates modifying the constraints in (\ref{gen-eq}) as:
\begin{equation}\label{gen-iq}
\tag{IQ-G}
\begin{aligned}
    \mathbf{z}_{1} &\geq \phi_{1}(\mathbf{x}),\\
    \mathbf{z}_{i+1} &\geq \phi_{i+1}(\mathbf{x},\boldsymbol{\omega}_i), \ i=1,\dots,L-2.
\end{aligned}
\end{equation}

\begin{proposition}\label{prob-cvx}
Given $\mathbf{x}$, we define the sets $E(\mathbf{x}):=\{\boldsymbol{\omega}_{L-1}|\boldsymbol{\omega}_{L-1} \text{ satisfies }(\ref{gen-eq})\},\\I(x):=\{\boldsymbol{\omega}_{L-1}|\omega_{L-1} \text{ satisfies }(\ref{gen-iq})\}$. Under Assumption 1, $R_{\boldsymbol{\theta}}$ defined by (\ref{gen-eq}) satisfies 
\begin{equation}
    R_{\boldsymbol{\theta}}(\mathbf{x})=\inf_{\boldsymbol{\omega}_{L-1}\in E(\mathbf{x})}\phi_L(\mathbf{x},\boldsymbol{\omega}_{L-1})=\inf_{\boldsymbol{\omega}_{L-1}\in I(\mathbf{x})}\phi_L(\mathbf{x},\boldsymbol{\omega}_{L-1}).
\end{equation}
\end{proposition}

\begin{proof}
Note that $E(\mathbf{x})$ is a singleton that consists of $\hat{\boldsymbol{\omega}}_{L-1}$ which satisfy (\ref{gen-eq}) given $\mathbf{x}$. Hence, $R(\mathbf{x})=\inf_{\omega_{L-1}\in E(\mathbf{x})}\phi_L(\mathbf{x},\boldsymbol{\omega}_{L-1})$.
Since $E(\mathbf{x})\subset I(\mathbf{x})$, so  $\inf_{\boldsymbol{\omega}_{L-1}\in E(\mathbf{x})}\phi_L(x,\boldsymbol{\omega}_{L-1})\geq\inf_{\boldsymbol{\omega}_{L-1}\in I(\mathbf{x})}\phi_L(\mathbf{x},\boldsymbol{\omega}_{L-1})$.\\
For $\boldsymbol{\omega}_{L-1}=(\mathbf{z}_1,\dots,\mathbf{z}_{L-1})\in I(\mathbf{x})$, we have $\mathbf{z}_1\geq\phi_1(\mathbf{x})=\hat{\mathbf{z}}_1$, $\mathbf{z}_i\geq\phi_i(\mathbf{x},\boldsymbol{\omega}_i)=\hat{\mathbf{z}}_i$ for $i=2,\dots,L-1$, where $\hat{\boldsymbol{\omega}}_{L-1}=(\hat{\mathbf{z}}_1,\dots,\hat{\mathbf{z}}_{L-1})$. Therefore, $\hat{\boldsymbol{\omega}}_{L-1}\leq \boldsymbol{\omega}_{L-1}$ for all $\boldsymbol{\omega}_{L-1}\in I(\mathbf{x})$. Since $\phi_L^\mathbf{x}$ is non-decreasing, we have:
\[\phi_L(\mathbf{x},\hat{\boldsymbol{\omega}}_{L-1})\leq \phi_L(\mathbf{x},\boldsymbol{\omega}_{L-1}).\]
Therefore, $\inf_{\boldsymbol{\omega}_{L-1}\in E(\mathbf{x})}\phi_L(\mathbf{x},\boldsymbol{\omega}_{L-1})\leq\inf_{\boldsymbol{\omega}_{L-1}\in I(\mathbf{x})}\phi_L(\mathbf{x},\boldsymbol{\omega}_{L-1})$. Combining with the other inequality, this shows that $\inf_{\boldsymbol{\omega}_{L-1}\in E(\mathbf{x})}\phi_L(\mathbf{x},\boldsymbol{\omega}_{L-1})=\inf_{\boldsymbol{\omega}_{L-1}\in I(\mathbf{x})}\phi_L(\mathbf{x},\boldsymbol{\omega}_{L-1})$. 
\end{proof}

We make the following assumption on the data fidelity and the regularization parameter.

\textbf{Assumption 2.} 
$D$ is convex in $\mathbf{x}$ and $\gamma>0$.

\begin{theorem}\label{ref-eq}
Under Assumptions 1 and 2, the following problem is convex
\begin{equation}\label{var-gnew}
    \min_{\mathbf{x},\boldsymbol{\omega}_{L-1}}D(\mathbf{Ax,y})+\gamma\phi_L(\mathbf{x},\boldsymbol{\omega}_{L-1}) \text{  subject to $\boldsymbol{\omega}_{L-1}$ satisfies (\ref{gen-iq})}.
\end{equation}
Furthermore, we denote $\mathcal{S}_1$ as set of minimizers of (\ref{var-prob}) with $R_{\boldsymbol{\theta}}$ defined by (\ref{gen-eq}), and $\mathcal{S}_2$ as set of minimizers of (\ref{var-gnew}). Then $\mathbf{x}\in \mathcal{S}_1$ if and only if there exists $\boldsymbol{\omega}_{L-1}$ such that $(\mathbf{x},\boldsymbol{\omega}_{L-1})\in\mathcal{S}_2$.
\end{theorem}

\begin{proof}
The constraints (\ref{gen-iq}) are convex since $\phi_i$ are convex. In particular, $D$ and $\gamma\phi_L$ are convex in $\mathbf{x}$, then problem (\ref{var-gnew}) is convex.
Due to proposition \ref{prob-cvx}, we have:
\begin{equation}
    \begin{aligned}
&\min_\mathbf{x}D(\mathbf{Ax,y})+\gamma R_{\boldsymbol{\theta}}(\mathbf{x})=\min_\mathbf{x}D(\mathbf{Ax,y})+\gamma\inf_{\boldsymbol{\omega}_{L-1}\in I(\mathbf{x})}\phi_L(\mathbf{x},\boldsymbol{\omega}_{L-1})\\&=\min_{\mathbf{x},\boldsymbol{\omega}_{L-1}\in I(\mathbf{x})}D(\mathbf{Ax,y})+\gamma\phi(\mathbf{x},\boldsymbol{\omega}_{L-1}).
    \end{aligned}
\end{equation}

Therefore, we have $\mathbf{x}\in\mathcal{S}_1$ if and only if there exists $\boldsymbol{\omega}_{L-1}$ such that $(\mathbf{x},\boldsymbol{\omega}_{L-1})\in\mathcal{S}_2$.
\end{proof}

\begin{corollary}\label{var-eq}
Consider the problem (\ref{var-prob}) with $R_{\boldsymbol{\theta}}$ given by an ICNN. Under Assumption 2, the following problem is convex
\begin{equation}\label{var-2}
\tag{P1}
\begin{aligned}
    \min_{\mathbf{x,z}_1,\dots,\mathbf{z}_{L-1}}&D(\mathbf{Ax,y})+\gamma h_L(\mathbf{V}_{L-1}\mathbf{x}+\mathbf{W}_{L-1}\mathbf{z}_{L-1}+\mathbf{b}_{L-1})\\
    \text{subject to }
    &\mathbf{z}_{1} \geq h_{1}(\mathbf{V}_{0}\mathbf{x}+\mathbf{b}_0),\\
    &\mathbf{z}_{i+1}\geq h_{i+1}(\mathbf{V}_{i}\mathbf{x}+\mathbf{W}_{i}\mathbf{z}_i+\mathbf{b}_i), \ i=1,\dots,L-2.
\end{aligned}
\end{equation}
Furthermore, $\mathbf{x}$ is a minimizer of (\ref{var-prob}) if and only if there exists $\mathbf{z}_1,\mathbf{z}_2,\dots,\mathbf{z}_{L-1}$ such that $(\mathbf{x,z}_1,\mathbf{z}_2,\dots,\mathbf{z}_{L-1})$ is a minimizer of (\ref{var-2}).
\end{corollary}

\begin{proof}
We note that (\ref{icnn-arch}) is a special case of (\ref{gen-eq}) with $\phi_{i+1}(\mathbf{x},\boldsymbol{\omega}_i)=h_{i+1}(\mathbf{V}_i\mathbf{x}+\mathbf{W}_i\mathbf{z}_i+\mathbf{b}_i)$. Therefore, Assumption 1 is satisfied in the ICNN setting. The result directly follows from Theorem \ref{ref-eq} with $(\mathbf{z}_1,\dots,\mathbf{z}_{L-1})=:\boldsymbol{\omega}_{L-1}$. 
\end{proof}

\section{Primal-Dual framework}
\label{pd-framework}
The final step of the proposed framework for solving (\ref{var-2}) is to replace the inequality constraints by indicator functions and reformulate (\ref{var-2}) as an equivalent unconstrained problem: 
\begin{equation}\label{var-3}
    \begin{aligned}
        \min_{\mathbf{x,z}_1,\dots,\mathbf{z}_{L-1}}D(\mathbf{Ax,y})+\gamma h_L(\mathbf{V}_{L-1}\mathbf{x}+\mathbf{W}_{L-1}\mathbf{z}_{L-1}+\mathbf{b}_{L-1})\\+\delta_{C_1}(\mathbf{V}_{0}\mathbf{x}+\mathbf{b}_0,\mathbf{z}_1)+\sum_{i=2}^{L-1}\delta_{C_i}(\mathbf{V}_{i-1}\mathbf{x}+\mathbf{W}_{i-1}\mathbf{z}_{i-1}+\mathbf{b}_{i-1},\mathbf{z}_{i}),
    \end{aligned}
\end{equation}
here $C_i:=\{(p,q)|h_{i}(p)\leq q\}$, and the indicator function is given by $\delta_{C_i}(\mathbf{x})$ which is $0$ if $(p,q)\in C_i$ and $\infty$ otherwise. We then apply a primal-dual algorithm to solve (\ref{var-3}).

To utilize the PDHG algorithm \cite{chambolle2011first,esser2010general,zhu2008efficient}, we recast (\ref{var-3}) in the following form:
\begin{equation}\label{pdhg-p}
    \min_\mathbf{u}\left\{\sum_{i=0}^Lf_i(\mathbf{K}_i\mathbf{u})+g(\mathbf{u})\right\}.
\end{equation}

We introduce the variable $\mathbf{u}=(\mathbf{x,z}_1,\dots,\mathbf{z}_{L-1})$ and consider:
\begin{equation}
    \begin{aligned}
        &\mathbf{K}_1=\left(\begin{array}{ccccc}
             \mathbf{V}_{0} & \mathbf{0} & \mathbf{0} & \cdots & \mathbf{0}\\
             \mathbf{0} & \mathbf{I} & \mathbf{0} & \cdots & \mathbf{0} 
        \end{array}\right)\\
        &\mathbf{K}_i=\left(\begin{array}{ccccccccc}
             \mathbf{V}_{i-1} & \mathbf{0} & \cdots & \mathbf{0} & \mathbf{W}_{i-1} & \mathbf{0} & \cdots & \mathbf{0} & \mathbf{0}\\
             \mathbf{0} & \mathbf{0} & \cdots & \mathbf{0} & \mathbf{0} & \mathbf{I} & \mathbf{0} & \cdots & \mathbf{0}  
        \end{array}\right), i=2,\dots,L-1\\
        &\mathbf{K}_L=\left(\begin{array}{ccccccccc}
             \mathbf{V}_{L-1} & \mathbf{0} & \cdots & \mathbf{0} & \mathbf{0} & \mathbf{0} & \cdots & \mathbf{0} & \mathbf{W}_{L-1}
        \end{array}\right)\\
        &\boldsymbol{\beta}_i=\left(\begin{array}{c}
            b_{i-1} \\
            \mathbf{0} \\
        \end{array}\right),i=1,\dots,L.
    \end{aligned}
\end{equation}
The data fidelity term $D(\mathbf{Ax,y})$ can either be included as $g(\mathbf{u})$ or as $f_0(\mathbf{K}_0\mathbf{u})$, where $\mathbf{K}_0=\left(\begin{array}{cccc}
\mathbf{A} & \mathbf{0} & \cdots & \mathbf{0} 
\end{array}\right)$. We also denote $\mathbf{K}^*_\mathbf{x},\mathbf{K}^*_{z_i}$ the rows of $\mathbf{K}^*$ with respect to the corresponding primal variables.

We then consider the following updates of PDHG:
\begin{equation}\label{pdhg}
        \begin{aligned}
            \mathbf{u}^{k+1} &= \text{prox}_g^\mathbf{T}(\mathbf{u}^k-\mathbf{TK}^*\mathbf{v}^k)\\
            \overline{\mathbf{u}}^{k+1} &= \mathbf{u}^{k+1}+\theta\left(\mathbf{u}^{k+1}-\mathbf{u}^k\right)\\
            \mathbf{v}^{k+1}_i &= \text{prox}_{f_i^*}^{\mathbf{S}_i}(\mathbf{v}_i^k+\mathbf{S}_i\mathbf{K}_i\overline{\mathbf{u}}^{k+1}),\ i=0,\dots,L, 
        \end{aligned}
\end{equation}
here the proximal operators are defined as $\text{prox}^\mathbf{S}_h(\mathbf{x})=\arg\min_{\mathbf{x}'}\left\{\frac{1}{2}\|\mathbf{x}'-\mathbf{x}\|^2_{\mathbf{S}^{-1}}+h(\mathbf{x})\right\}$, where $\|\mathbf{x}\|^2_{\mathbf{S}^{-1}}=\langle \mathbf{x},\mathbf{S}^{-1}\mathbf{x}\rangle$, and the step-size matrices $\mathbf{T,S}_i$ are symmetric and positive definite. The algorithm is known to converge \cite{pock2011diagonal} if $\|\mathbf{S}^{1/2}\mathbf{KT}^{1/2}\|<1$ and $\theta=1$, where $\mathbf{S}=\operatorname{diag}(\mathbf{S}_0,\dots,\mathbf{S}_L)$. We choose diagonal matrices $\mathbf{T}=\operatorname{diag}(\tau_0,\dots,\tau_{L-1}),\ \mathbf{S}_i=\operatorname{diag}(\sigma_i\mathbf{I},\sigma_i\mathbf{I})$ for $i=1,\dots,L-1$ and $\mathbf{S}_j=\sigma_j\mathbf{I}$ for $j=0,L$ as our step-size matrices. 

\rv{Note that in our setting, the function $g$ only depends on $\mathbf{x}$, hence the primal updates for variables $\mathbf{z}_i$ are simply affine transformations.} Moreover, the Moreau identity relates the proximal operator of a function $h$ to that of its conjugate $h^*$ defined by $h^*(\mathbf{y})=\sup_\mathbf{x}\langle \mathbf{x,y}\rangle-h(\mathbf{x})$, updates for $\mathbf{v}_i$s can be computed via $\text{prox}_{f_i}$, which are the projections onto $C_i$ or $\text{prox}_{h_L}$. With common choices of activations such as ReLU, leaky ReLU, these operators can be computed exactly.

\rv{Specifically, the dual updates for $i=1,\dots,L-1$ at entries $(p,q)$ can be computed with:}
\rv{
\begin{equation}
\operatorname{prox}_{f_i^*}^{\sigma_i}((\bar{p},\bar{q}))    =(\bar{p},\bar{q})-\sigma_i\left(\operatorname{prox}_{\delta_{C_i}}^{\sigma_i^{-1}}\left(\left(\frac{\bar{p}}{\sigma_i}+b_{i-1},\frac{\bar{q}}{\sigma_i}\right)\right)-b_{i-1}\right),
\end{equation}}
\rv{where we have applied the translation property of proximal operators and the Moreau identity. Consider leaky ReLU activation for example, the proximal operator of $\delta_{C_i}$ corresponds to epigraphical projection, which is given by:} 
\rv{
\begin{equation}
    \operatorname{proj}_{C_1}(\bar{p},\bar{q})=\begin{cases}
    (\bar{p},\bar{q}) & \text{if } h_1(\bar{p})\leq \bar{q}\\
    (\frac{\bar{p}+\bar{q}}{2},\frac{\bar{p}+\bar{q}}{2}) & \text{if }|\bar{q}|\leq \bar{p}\\
    (\frac{\bar{p}+\alpha\bar{q}}{1+\alpha^2},\frac{\alpha(\bar{p}+\alpha \bar{q})}{1+\alpha^2}) & \text{if }\bar{q}\leq\alpha \bar{p}\text{ and } \bar{p}\leq-\alpha \bar{q}\\
    (0,0) & \text{otherwise}
    \end{cases}.
\end{equation}
}
\rv{where $h_1$ denotes the leaky ReLU function with negative slope $\alpha$.
}

\rv{The general steps for the proposed method is then outlined in Algorithm \ref{alg:algo}.}

\begin{algorithm*}[h]
\rv{
\caption{Proposed method}
\label{alg:algo}
\begin{algorithmic}[1]
\State{Input: Step size matrices $\mathbf{T,S_i}$ for $i=0,\dots,L$, number of maximum iterations $N_{maxiter}$\\ Initialization for primal variables $\mathbf{u}^0=(\mathbf{x}^0,\mathbf{z}^0_1\dots,\mathbf{z}^0_{L-1})$ and dual variables $\mathbf{v}^0_i=(\mathbf{v}^0_{i1},\mathbf{v}^0_{i2})$ for $i=0,\dots,L-1$ and $\mathbf{v}^0_L$.}
\For{$k=0,1, \cdots N_{maxiter}-1$:}
\State{$\mathbf{x}^{k+1}=\operatorname{prox}^{\tau_0}_g(\mathbf{x}^k-\tau_0\mathbf{K}^*_\mathbf{x}\mathbf{v}^k)$}
\State{$\mathbf{z}^{k+1}_i=\mathbf{z}^k_i-\tau_i\mathbf{K}^*_{\mathbf{z}_i}\mathbf{v}^k$ for $i=1,\dots,L-1$}
\State{$\bar{\mathbf{u}}^{k+1}=(2\mathbf{x}^{k+1}-\mathbf{x}^{k},2\mathbf{z}_1^{k+1}-\mathbf{z}_1^{k}\dots,2\mathbf{z}_{L-1}^{k+1}-\mathbf{z}_{L-1}^{k})$}
\State{$\mathbf{v}^{k+1}_0=\operatorname{prox}^{\sigma_0}_{f^*_0}(\mathbf{v}^{k}_0+\sigma_0\mathbf{K}_0\bar{\mathbf{u}}^{k+1})$}
\State{$(\mathbf{v}^{k+1}_{i1},\mathbf{v}^{k+1}_{i2})=(\mathbf{v}^{k}_{i1},\mathbf{v}^{k}_{i2})-\sigma_i\left(\operatorname{prox}_{\delta_{C_i}}^{\sigma_i^{-1}}\left(\left(\frac{\mathbf{v}^{k}_{i1}}{\sigma_i}+b_{i-1},\frac{\mathbf{v}^{k}_{i2}}{\sigma_i}\right)\right)-b_{i-1}\right)$ for $i=1,\dots,L-1$}
\State{$\mathbf{v}^{k+1}_L=\operatorname{prox}^{\sigma_L}_{\gamma h_L}(\mathbf{v}^{k}_L+\sigma_L\mathbf{K}_L\bar{\mathbf{u}}^{k+1})$}
\EndFor
\State{Output: Reconstructed image $\mathbf{x}^{N_{maxiter}}$}
\end{algorithmic}
}
\end{algorithm*}

We note that the proposed primal-dual framework introduces auxiliary variables $\mathbf{z}_i$. However, these auxiliary variables correspond directly to the layer-wise activations already present in the network. Hence, the method does not incur additional memory costs compared to standard backpropagation \cite{li2019lifted}. Moreover, the updates for the auxiliary variables and the dual variables can be computed independently, which offers the potential for efficient parallel computation.

\section[Experiments]{Experiments \footnote{The python code are available from https://doi.org/10.5281/zenodo.17426033.}}
\label{sec:experiments}
We evaluate the performance of the proposed method and compare with subgradient methods on three imaging tasks, (i) salt and pepper denoising, (ii) sparse-view CT reconstruction, and (iii) image inpainting. For all tasks, we utilize a learned regularizer parametrized by an ICNN, which consists of a convolution layer and a global averaging operator layer, followed by two fully connected layers. The regularizer can be represented by $R_{\boldsymbol{\theta}}(\mathbf{x})=\mathbf{W}_2h_2(\mathbf{W}_1\mathbf{P}\mathbf{z}+\mathbf{b}_1)$ with $\mathbf{z}=h_1(\mathbf{V}_0x+\mathbf{b}_0)$. Here $\mathbf{V}_0$ corresponds to a convolution operator with 32 $5\times5$ filters, and $\mathbf{P}$ denotes an average pooling operater with $16\times 16$ pool size. The fully connected layers $\mathbf{W}_1,\mathbf{W}_2$ consists of 256 and one output neurons respectively. The activations $h_1,h_2$ are chosen to be leaky ReLU and ReLU respectively, with the leaky ReLU's negative slope set to $0.2$. \rv{Variants of this baseline architecture with additional layers or smooth activations are also considered in the CT reconstruction and inpainting experiments, respectively.}

We follow the adversarial training framework in \cite{lunz2018adversarial,mukherjee2020learned} to train the regularizer $R_{\boldsymbol{\theta}}$. \rv{For non-smooth regularizers, we compute $\nabla R_{\boldsymbol{\theta}}$ in (\ref{eq:AR_loss}) as a subgradient via automatic differentiation. The training objective (\ref{eq:AR_loss}) is optimized with the Adam optimizer \cite{adam2014method}. The gradient penalty parameter $\lambda_{\mathrm{GP}}$ is selected empirically to optimize reconstruction quality as measured by PSNR. The specific training parameters vary across the three tasks and are detailed in the corresponding subsections.} 

The associated minimization problem is then solved with the proposed method, and compare with the subgradient method with (a) constant step-size (SM-C) and (b) diminishing step-size (SM-D), with step-size at the $k$-th iteration given by the initial step-size divided by $k$. The corresponding subgrdradient methods SM-C and SM-D are given in Algorithm \ref{alg:smc} and \ref{alg:smd} respectively. For both methods, the subgradients are computed using automatic differentiation. 
\begin{algorithm*}[h!]
\caption{Subgradient method with constant step-size (SM-C)  \cite{boyd2003subgradient}}
\label{alg:smc}
\begin{algorithmic}[1]
\State{Input: Initialization $x^0$, constant step-size $\eta$, maximum number of iterations $N_{max}$.}
\For{$k=0,1,\dots,N_{max}$}
\State{Compute subgradient $\mathbf{g}^k\in\partial_{\mathbf{x}}(D(\mathbf{Ax}^k,\mathbf{y})+\gamma R_{\boldsymbol{\theta}}(\mathbf{x}^k))$.}
\State{Update $\mathbf{x}^{k+1}=\mathbf{x}^k-\eta\mathbf{g}^k$.}
\EndFor
\end{algorithmic}
\end{algorithm*}

\begin{algorithm*}[h!]
\caption{Subgradient method with diminishing step-size (SM-D) \cite{boyd2003subgradient}}
\label{alg:smd}
\begin{algorithmic}[1]
\State{Input: Initialization $x^0$, initial step-size $\eta^0$, maximum number of iterations $N_{max}$.}
\For{$k=0, 1,\dots,N_{max}$}
\State{Compute subgradient $\mathbf{g}^k\in\partial_{\mathbf{x}}(D(\mathbf{Ax}^k,\mathbf{y})+\gamma R_{\boldsymbol{\theta}}(\mathbf{x}^k))$.}
\State{Update $\mathbf{x}^{k+1}=\mathbf{x}^k-\eta^k\mathbf{g}^k$.}
\State{Update step-size $\eta^{k+1}=\eta^k/k$.}
\EndFor
\end{algorithmic}
\end{algorithm*}

\subsection{Salt and pepper denoising}
In this example, 1000 grayscale images from the FFHQ dataset \cite{karras2019style} downsampled to size $256\times256$ are used as training \rv{ground truth data}. The salt and pepper corrupted images are used as degraded samples, i.e. $\mathbf{A}^\dagger\mathbf{y}_j$ for adversarial training. \rv{The regularizer is trained for 20 epochs with a learning rate of $5e$-$4$, $(\beta_1,\beta_2)=(0.5,0.99)$, batch size of 8, and $\lambda_{GP}=5$.} To deal with salt and pepper noise, we utilize an $L^1$-data fidelity \cite{chambolle2011first}. The optimization problem is formulated as:
\begin{equation}
    \min_{\mathbf{x,z}}\lambda\|\mathbf{x}-\mathbf{y}\|_1+\mathbf{W}_2h_2(\mathbf{W}_1\mathbf{P}\mathbf{z}+\mathbf{b}_1)+\delta_{C_1}(\mathbf{V}_0\mathbf{x}+\mathbf{b}_0,\mathbf{z}),
\end{equation}
where $C_1=\{(p,q)|h_1(p)\leq q\}$.
The steps to solve the variational problem are outlined as follows:

\begin{equation}\label{pdhg-l1}
        \begin{aligned}
            \mathbf{x}^{k+1},\mathbf{z}^{k+1} &= \text{prox}_{\lambda\|\cdot-\mathbf{y}\|_1}^{\tau_1}(\mathbf{x}^k-\tau_1\mathbf{V}_0^*\mathbf{v}^k_{1,1}),\mathbf{z}^k-\tau_2(\mathbf{v}^k_{1,2}+\mathbf{P}^*\mathbf{W}_1^*\mathbf{v}^k_2)\\\overline{\mathbf{x}}^{k+1},\overline{\mathbf{z}}^{k+1} &= 2\mathbf{x}^{k+1}-\mathbf{x}^k,2\mathbf{z}^{k+1}-\mathbf{z}^k\\(\tilde{\mathbf{v}}^{k+1}_{1,1},\tilde{\mathbf{v}}^{k+1}_{1,2}),\tilde{\mathbf{v}}^{k+1}_2 &=(\mathbf{v}^k_{1,1}+\sigma_1\mathbf{V}_0\overline{\mathbf{x}}^{k+1},\mathbf{v}^k_{1,2}+\sigma_1\overline{\mathbf{z}}^{k+1}),\mathbf{v}^k_2+\sigma_2\mathbf{W}_1\mathbf{P}\overline{\mathbf{z}}^{k+1}\\
            \mathbf{v}^{k+1}_1,\mathbf{v}^{k+1}_2 &= \tilde{\mathbf{v}}^{k+1}_1-\sigma\text{proj}_{C_1}\left(\frac{\tilde{\mathbf{v}}^{k+1}_1}{\sigma}+\boldsymbol{\beta}_1\right),\tilde{\mathbf{v}}^{k+1}_2- \text{prox}_{f_2}^{\sigma_2^{-1}}\left(\frac{\tilde{\mathbf{v}}^{k+1}_2}{\sigma_2}\right),
        \end{aligned}
\end{equation}
where the proximal operator of the $L^1$ data fidelity is given by the pointwise soft shrinkage function:
\begin{equation}
\left[\operatorname{prox}_{\lambda\|\cdot-\mathbf{y}\|_2}^{\tau}(\bar{\mathbf{x}})\right]_i=\begin{cases}
        \bar{\mathbf{x}}_i-\tau\lambda &\text{if }\bar{\mathbf{x}}_i-\mathbf{y}_i>\tau\lambda\\
        \bar{\mathbf{x}}_i+\tau\lambda &\text{if }\bar{\mathbf{x}}_i-\mathbf{y}_i<-\tau\lambda\\
        \mathbf{y}_i &\text{otherwise}
    \end{cases}.
\end{equation}

We consider vector-valued step-sizes, $\mathbf{T}=\operatorname{diag}(\tau_1\mathbf{I}_\mathbf{x},\tau_2\mathbf{I}_\mathbf{z}),S=\operatorname{diag}(\sigma_1\mathbf{I}_{\mathbf{v}_1},\sigma_2\mathbf{I}_{\mathbf{v}_2})$. In terms of the general form, we incorporate the data fidelity as $g$, and consider the operator:
\begin{equation*}
    \mathbf{K}=\left(\begin{array}{cc}
             \mathbf{V}_0 & 0\\
             0 & \mathbf{I}\\
             0 & \mathbf{W}_1\mathbf{P}
        \end{array}\right).
\end{equation*}
Then $\mathbf{S}^{1/2}\mathbf{KT}^{1/2}$ is given by 
\begin{equation*}
    \mathbf{K}=\left(\begin{array}{cc}
             \sqrt{\sigma_1\tau_1}\mathbf{V}_0 & 0\\
             0 & \sqrt{\sigma_1\tau_2}\mathbf{I}\\
             0 & \sqrt{\sigma_2\tau_2}\mathbf{W}_1\mathbf{P}
        \end{array}\right).
\end{equation*}
To study the convergence condition we compute 
\begin{equation*}
    \begin{aligned}
        \|\mathbf{S}^{1/2}\mathbf{KT}^{1/2}\mathbf{u}\|^2&=\sigma_1\tau_1\|\mathbf{V}_0\mathbf{x}\|^2+\sigma_1\tau_2\|\mathbf{z}\|^2+\sigma_2\tau_2\|\mathbf{W}_1\mathbf{Pz}\|^2\\
        &\leq \sigma_1\tau_1\|\mathbf{V}_0\|^2\|\mathbf{x}\|^2+\sigma_1\tau_2\|\mathbf{z}\|^2+\sigma_2\tau_2\|\mathbf{W}_1\mathbf{P}\|^2\|\mathbf{z}\|^2\\
        &=\sigma_1\tau_1\|\mathbf{V}_0\|^2\|\mathbf{x}\|^2+(\sigma_1\tau_2+\sigma_2\tau_2\|\mathbf{W}_1\mathbf{P}\|^2)\|\mathbf{z}\|^2.
    \end{aligned}
\end{equation*}
The step-sizes are then chosen based on the condition $\|\mathbf{S}^{1/2}\mathbf{KT}^{1/2}\|<1$ and are given by:
\begin{equation}\label{step-size}
    \sigma_1=\frac{c_1}{\|\mathbf{V}_0\|^2},\sigma_2=\frac{c_2}{\|\mathbf{W}_1\mathbf{P}\|^2},\tau_1=\frac{1}{\sigma_1\|\mathbf{V}_0\|^2},\tau_2=\frac{1}{\sigma_1+\sigma_2\|\mathbf{W}_1\mathbf{P}\|^2},
\end{equation}
with hyperparameters $c_1,c_2$.

\begin{figure}[!h]
   \begin{minipage}[b]{0.45\textwidth}
    \centering     
         \begin{tikzpicture}
    \node[anchor=south west,inner sep=0] (image) at (0,0) {\includegraphics[width=\textwidth]{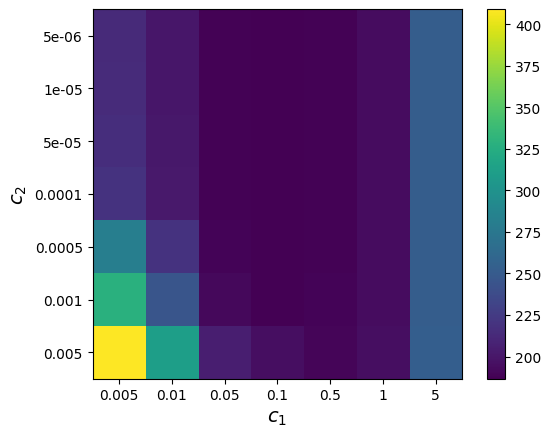}};
    \begin{scope}[x={(image.south east)},y={(image.north west)}]
    \draw[line width=.1mm,yellow] (0.215,0.67) node[regular polygon, regular polygon 
    sides=4,fill=yellow,scale=0.75,draw] {};
    \draw[line width=.1mm,purple] (0.505,0.67) node[star,fill=purple,scale=0.75,draw] {};
    \draw[line width=.01mm,blue] (0.505,0.91) node[regular polygon, regular polygon sides=6,fill=blue,scale=0.75,draw] {};
    \draw[line width=.1mm,green] (0.795,0.67) node[regular polygon, regular polygon sides=5,fill=green,scale=0.75,draw] {};;
    \draw[line width=.1mm,red] (0.505,0.19) node[diamond,fill=red,scale=0.75,draw] {};
    \end{scope}
\end{tikzpicture}
     \centerline{Average objective value}
   \end{minipage}
   \hspace{1cm}
   \begin{minipage}[b]{0.45\textwidth}
     \centering
     \includegraphics[width=\linewidth]{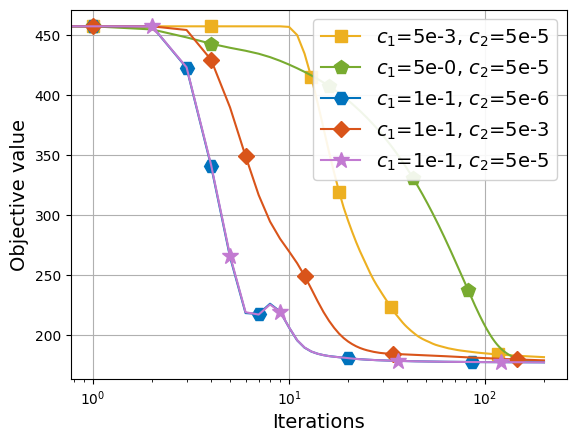}
     \centerline{Objective value}
   \end{minipage}      
   \caption{Denoising: Ablation study of proposed method for step-size hyperparameters. The markers on the left corresponds to those depicted in the energy versus iterations plots on the right.}
   \label{fig:sp_tune}
\end{figure}
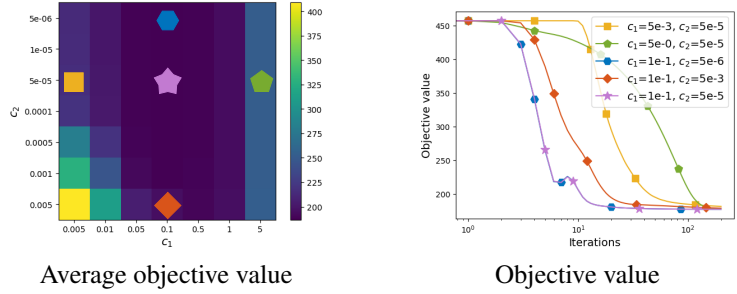

\begin{figure}[!h]
\begin{minipage}{0.45\textwidth}
     \centering
     \includegraphics[width=\linewidth]{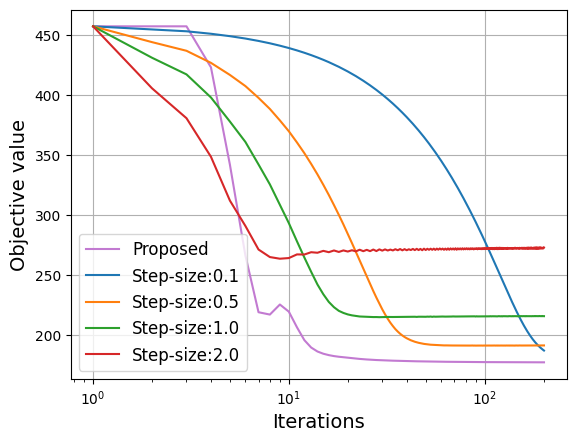}
   \end{minipage}
   \hspace{1cm}
    \begin{minipage}{0.45\textwidth}
     \centering
     \includegraphics[width=\linewidth]{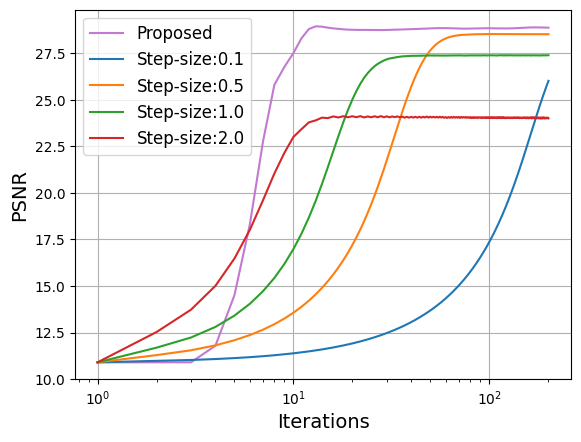}
   \end{minipage}    

\begin{minipage}{0.45\textwidth}
    \centering
     \includegraphics[width=\linewidth]{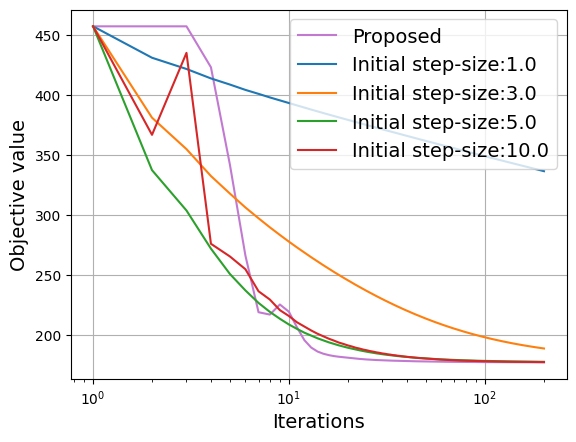}
     \hspace{1cm}
     \centerline{SM-C}
   \end{minipage}
   \hspace{1cm}
      \begin{minipage}{0.45\textwidth}
     \centering
     \includegraphics[width=\linewidth]{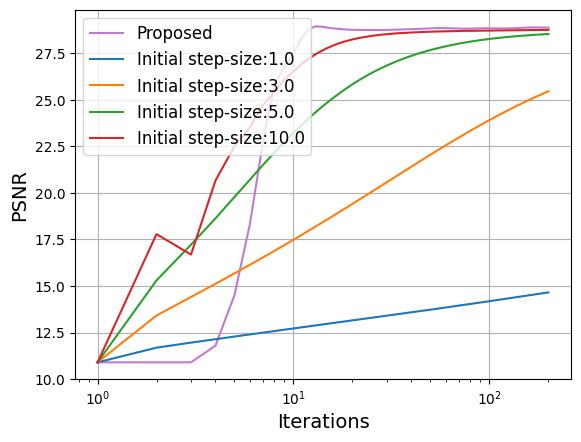}
     \centerline{SM-D}
   \end{minipage}    
\caption{Denoising: Comparison to subgradient methods. Note that SM-C either converges only to a suboptimal solution or is very slowly. SM-D converges, but remains slightly slower than the proposed method.}
   \label{fig:sp_comp}
\end{figure}

\textbf{Parameters:} For this experiment, we set $\lambda=0.02$. For the proposed method, we pick $c_1,c_2$ from \{$5e$-$3$,$1e$-$2$,$5e$-$2$,$1e$-$1$,$5e$-$1$,$1$,$5$\}, \{$5e$-$6$,$1e$-$5$,$5e$-$5$,$1e$-$4$,$5e$-$4$,$1e$-$3$,$5e$-$3$\}. For SM-C, we choose the step-size from $\{0.1,0.5,1,2\}$. As for SM-D, we select the initial step-size from $\{1,3,5,10\}$. 

\textbf{Ablation study: }Figure \ref{fig:sp_tune} shows the ablation study of the step-size hyperparameters for the proposed method. We ran 200 iterations of the proposed method for each hyperparameter combination and evaluated the average objective value to assess convergence. The left plot shows the average objective values, while the right plot depicts energy versus iterations for different values of $c_1,c_2$.

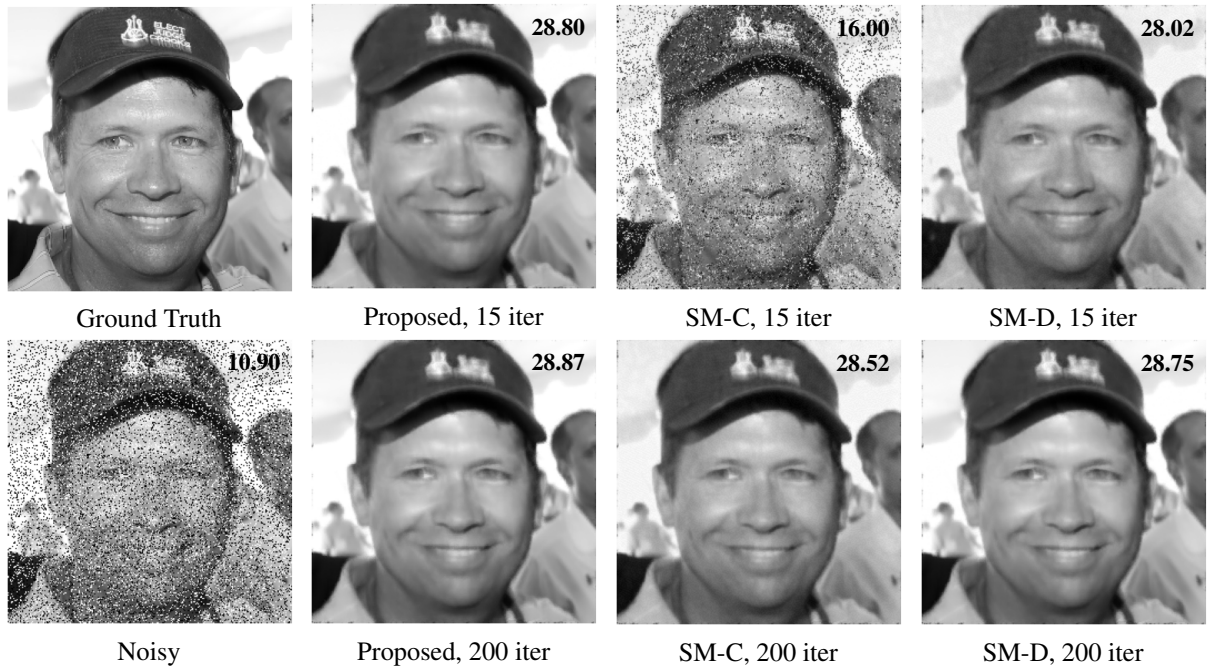
\begin{figure}[!h]
   \begin{minipage}{0.24\textwidth}
     \centering         
    \begin{tikzpicture}
    \node[anchor=south west,inner sep=0] (image) at (0,0) {\includegraphics[width=\textwidth]{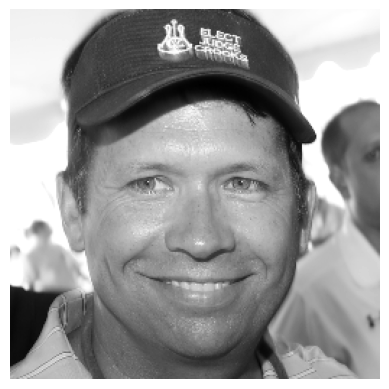}};
    \begin{scope}[x={(image.south east)},y={(image.north west)}]
    \node[text=black] at (0.85,0.9) {\footnotesize };
    \end{scope}
\end{tikzpicture}
     \centerline{Ground Truth}
   \end{minipage}
   \begin{minipage}{0.24\textwidth}
    \centering
    \begin{tikzpicture}
    \node[anchor=south west,inner sep=0] (image) at (0,0) {\includegraphics[width=\textwidth]{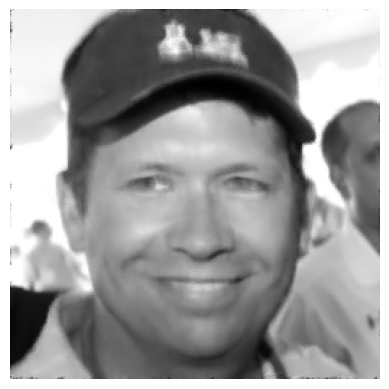}};
    \begin{scope}[x={(image.south east)},y={(image.north west)}]
    \node[text=black] at (0.85,0.9) {\footnotesize \textbf{28.80}};
    \end{scope}
\end{tikzpicture}
     \centerline{Proposed, 15 iter}
   \end{minipage}
   \begin{minipage}{0.24\textwidth}
    \centering
    \begin{tikzpicture}
    \node[anchor=south west,inner sep=0] (image) at (0,0) {\includegraphics[width=\textwidth]{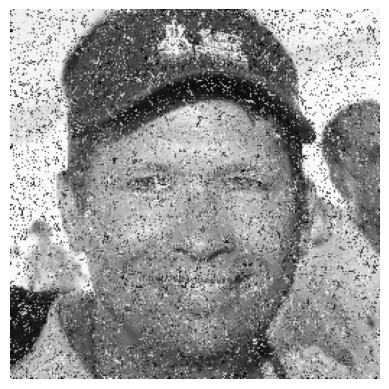}};
    \begin{scope}[x={(image.south east)},y={(image.north west)}]
    \node[text=black] at (0.85,0.9) {\footnotesize \textbf{16.00}};
    \end{scope}
\end{tikzpicture}
     \centerline{SM-C, 15 iter}
   \end{minipage}
   \begin{minipage}{0.24\textwidth}
    \centering
    \begin{tikzpicture}
    \node[anchor=south west,inner sep=0] (image) at (0,0) {\includegraphics[width=\textwidth]{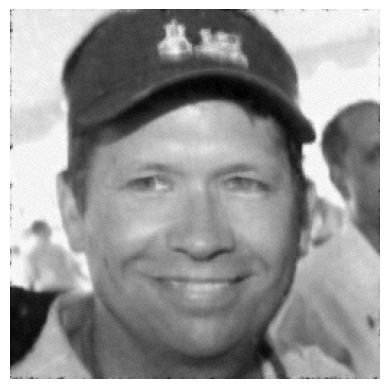}};
    \begin{scope}[x={(image.south east)},y={(image.north west)}]
    \node[text=black] at (0.85,0.9) {\footnotesize \textbf{28.02}};
    \end{scope}
\end{tikzpicture}
     \centerline{SM-D, 15 iter}
   \end{minipage}
   
    \begin{minipage}{0.24\textwidth}
     \centering
    \begin{tikzpicture}
    \node[anchor=south west,inner sep=0] (image) at (0,0) {\includegraphics[width=\textwidth]{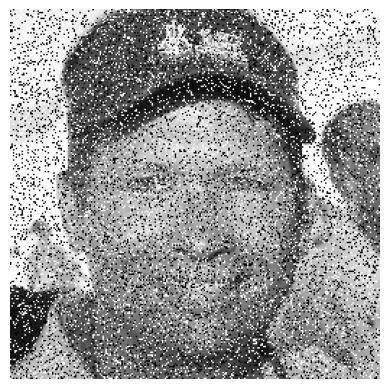}};
    \begin{scope}[x={(image.south east)},y={(image.north west)}]
    \node[text=black] at (0.85,0.9) {\footnotesize \textbf{10.90}};
    \end{scope}
\end{tikzpicture}
     \centerline{Noisy}
   \end{minipage}
      \begin{minipage}{0.24\textwidth}
     \centering
    \begin{tikzpicture}
    \node[anchor=south west,inner sep=0] (image) at (0,0) {\includegraphics[width=\textwidth]{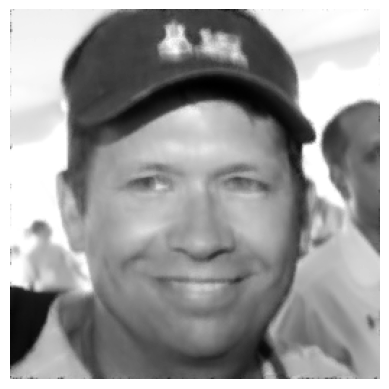}};
    \begin{scope}[x={(image.south east)},y={(image.north west)}]
    \node[text=black] at (0.85,0.9) {\footnotesize \textbf{28.87}};
    \end{scope}
\end{tikzpicture}
     \centerline{Proposed, 200 iter}
   \end{minipage}
      \begin{minipage}{0.24\textwidth}
     \centering
    \begin{tikzpicture}
    \node[anchor=south west,inner sep=0] (image) at (0,0) {\includegraphics[width=\textwidth]{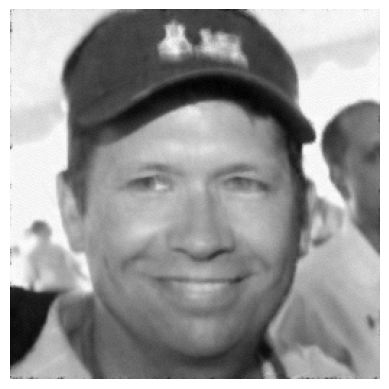}};
    \begin{scope}[x={(image.south east)},y={(image.north west)}]
    \node[text=black] at (0.85,0.9) {\footnotesize \textbf{28.52}};
    \end{scope}
\end{tikzpicture}
     \centerline{SM-C, 200 iter}
   \end{minipage}
      \begin{minipage}{0.24\textwidth}
     \centering
    \begin{tikzpicture}
    \node[anchor=south west,inner sep=0] (image) at (0,0) {\includegraphics[width=\textwidth]{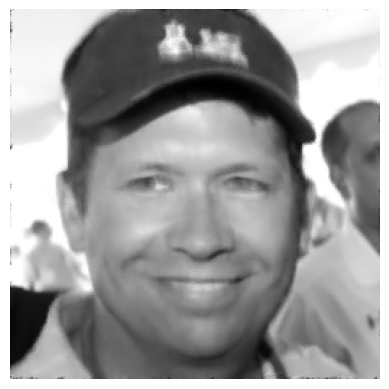}};
    \begin{scope}[x={(image.south east)},y={(image.north west)}]
    \node[text=black] at (0.85,0.9) {\footnotesize \textbf{28.75}};
    \end{scope}
\end{tikzpicture}
     \centerline{SM-D, 200 iter}
   \end{minipage}
   \caption{Denoising: Visual comparison of reconstructions, with PSNR values shown in the top right corner. The proposed method achieves a visually satisfactory reconstruction within 15 iterations, while that of SM-C remains noisy.}
   \label{fig:sp_recon}
\end{figure}

\textbf{Results:} The proposed method with the optimal choice of $c_1,c_2$ is then compared with the subgradient methods. The first row of Figure \ref{fig:sp_comp} shows energy versus iterations plots, indicating that the constant step-size subgradient methods fail to converge within 200 iterations, while the diminishing step-size subgradient methods do converge, albeit slower than the proposed method. Moreover, we evaluate the Peak Signal-to-Noise Ratio (PSNR). The proposed method achieves the highest PSNR values in less than 20 iterations, outperforming both subgradient methods. Figure \ref{fig:sp_recon} shows the reconstructed images produced by each method. Reconstructions are provided at both 15 and 200 iterations, with the proposed method delivering visually satisfactory results as early as 15 iterations. 

\subsection{CT with Poisson noise}
In this section, we consider a sparse-view computed tomography (CT) reconstruction task, with human abdominal CT scans of the Mayo clinic for the low-dose CT grand challenge \cite{mccollough2016tu} as training and testing data. The measurements are simulated using a parallel beam geometry with 200 angles and 400 bins. We model the noise as Poisson with a constant background level $r=50$, \rv{that is $\mathbf{y}\sim\text{Pois}(\mathbf{Ax}+r)$, where $\mathbf{A}$ is the scaled X-ray transform with the prescribed geometry. We take 2240 slices with batch size 20 from 9 patients as training ground truth data, and the filtered back projections (FBP) are used as degraded samples for training. The regularizer is trained for 20 epochs with a learning rate of of $2e$-$5$, $(\beta_1,\beta_2)=(0.5,0.99)$, batch size of 20, and $\lambda_{GP}=30$.} For reconstruction, we consider the  Kullback–Leibler (KL) divergence data fidelity, which is suitable for Poisson-distributed data: 
\begin{equation}
D(\mathbf{Ax,y})=\textbf{1}^T\left(\mathbf{Ax}-\mathbf{y}+r\right)+\mathbf{y}^T\log\left(\frac{\mathbf{y}}{\mathbf{Ax}+r}\right),
\end{equation}
where $\textbf{1}$ denotes a vector of 1s. We also assume that the reconstruction is bounded below by 0. The reformulation can then be written as:
\begin{equation}
    \min_{\mathbf{x,z}}D(\mathbf{Ax,y})+\gamma \mathbf{W}_2h_2(\mathbf{W}_1\mathbf{Pz}+\mathbf{b}_1)+\delta_{C_1}(\mathbf{V}_0\mathbf{x}+\mathbf{b}_0,\mathbf{z})+\delta_{[0,\infty)}(\mathbf{x}).
\end{equation}

Unlike the previous experiment, we dualize the forward operater $\mathbf{A}$ with the data fidelity acting as $f_0$. This leads to the following updates:
\begin{equation}\label{pdhg-ct}
        \begin{aligned}
            \mathbf{x}^{k+1},\mathbf{z}^{k+1}&=\text{max}(\mathbf{x}^k-\tau_1\mathbf{A}^*\mathbf{v}_0+\mathbf{V}_0^*\mathbf{v}^k_{1,1},0),\mathbf{z}^k-\tau_2(\mathbf{v}^k_{1,2}+\mathbf{P}^*\mathbf{W}_1^*\mathbf{v}^k_2)\\
            \overline{\mathbf{x}}^{k+1},\overline{\mathbf{z}}^{k+1}&=2\mathbf{x}^{k+1}-\mathbf{x}^k,2\mathbf{z}^{k+1}-\mathbf{z}^k\\            (\tilde{\mathbf{v}}^{k+1}_{1,1},\tilde{\mathbf{v}}^{k+1}_{1,2}),\tilde{\mathbf{v}}^{k+1}_2&=(\mathbf{v}^k_{1,1}+\sigma_1\mathbf{V}_0\overline{\mathbf{x}}^{k+1},\mathbf{v}^k_{1,2}+\sigma_1\overline{\mathbf{z}}^{k+1}),\mathbf{v}^k_2+\sigma_2\mathbf{W}_1\mathbf{P}\overline{\mathbf{z}}^{k+1}\\
            \mathbf{v}^{k+1}_0&=\text{prox}_{f_0^*}^{\sigma_0}(\mathbf{v}^k_0+\sigma_0\mathbf{A}\overline{\mathbf{x}}^{k+1})\\
            \mathbf{v}^{k+1}_1,\mathbf{v}^{k+1}_2 &= \tilde{\mathbf{v}}^{k+1}_1-\sigma\text{proj}_{C_1}\left(\frac{\tilde{\mathbf{v}}^{k+1}_1}{\sigma}+\beta_1\right),\tilde{\mathbf{v}}^{k+1}_2- \text{prox}_{f_2}^{\sigma_2^{-1}}\left(\frac{\tilde{\mathbf{v}}^{k+1}_2}{\sigma_2}\right),
        \end{aligned}
\end{equation}
where for the Kullback–Leibler divergence, we let $f_0(\mathbf{w})=\textbf{1}^T\left(\mathbf{w-y+r}\right)+\mathbf{y}^T\log\left(\frac{\mathbf{y}}{\mathbf{w+r}}\right)$. The proximal operator of its conjugate can be given by \cite{chambolle2018stochastic}:
\begin{equation}
\left[\operatorname{prox}_{f_0^*}^{\sigma_0}(\bar{\mathbf{w}})\right]_i=\frac{1}{2}\left(\bar{\mathbf{w}}_i+1+\sigma_0\mathbf{r}_i-\sqrt{(\bar{\mathbf{w}}_i-1+\sigma_0\mathbf{r}_i)^2+4\sigma_0\mathbf{y}_i}\right)
\end{equation}

Following the step-size selection scheme and incorporating the forward operator $\mathbf{A}$ into the operator $\mathbf{K}$, we pick step-sizes $\mathbf{T}=\operatorname{diag}(\tau_1\mathbf{I}_\mathbf{x},\tau_2\mathbf{I}_\mathbf{z}),\mathbf{S}=\operatorname{diag}(\sigma_0\mathbf{I}_{\mathbf{v}_0},\sigma_1\mathbf{I}_{\mathbf{v}_1},\sigma_2\mathbf{I}_{\mathbf{v}_2})$ given by:
\begin{equation}\label{step-size-ct}
    \sigma_0=\frac{c_0}{\|\mathbf{A}\|^2},\sigma_1=\frac{c_1}{\|\mathbf{V}_0\|^2},\sigma_2=\frac{c_2}{\|\mathbf{W}_1\mathbf{P}\|^2},\tau_1=\frac{1}{c_0+c_1},\tau_2=\frac{1}{\sigma_1+c_2}.
\end{equation} 

\begin{figure}[!h]
\begin{minipage}{0.45\textwidth}
     \centering
     \includegraphics[width=\linewidth]{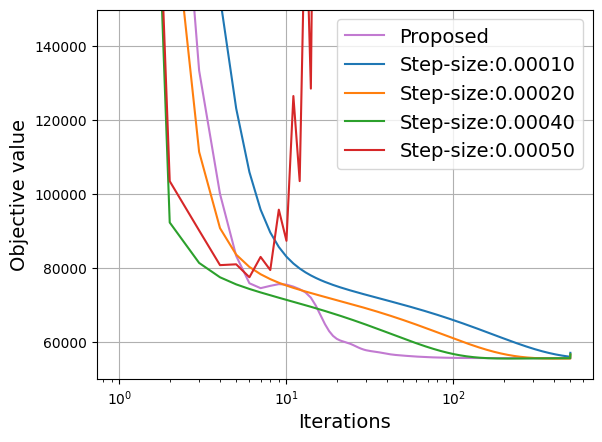}
   \end{minipage}
   \hspace{1cm}
    \begin{minipage}{0.45\textwidth}
     \centering
     \includegraphics[width=\linewidth]{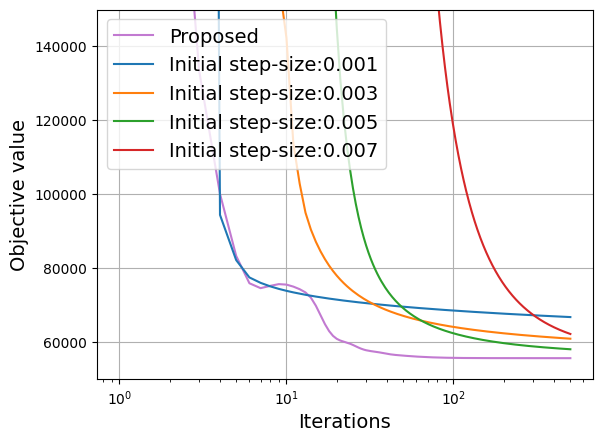}
   \end{minipage}    

\begin{minipage}{0.45\textwidth}
    \centering
     \includegraphics[width=\linewidth]{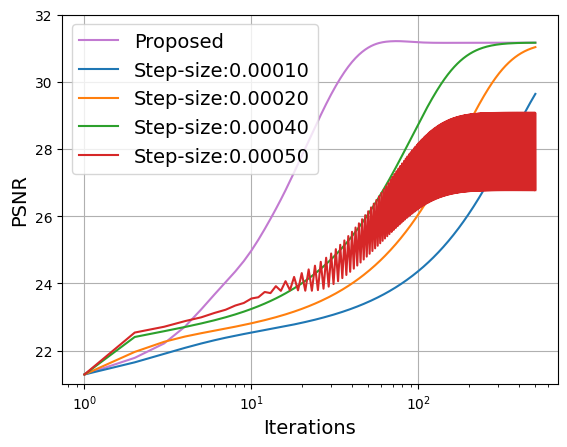}
     \hspace{1cm}
     \centerline{SM-C}
   \end{minipage}
   \hspace{1cm}
      \begin{minipage}{0.45\textwidth}
     \centering
     \includegraphics[width=\linewidth]{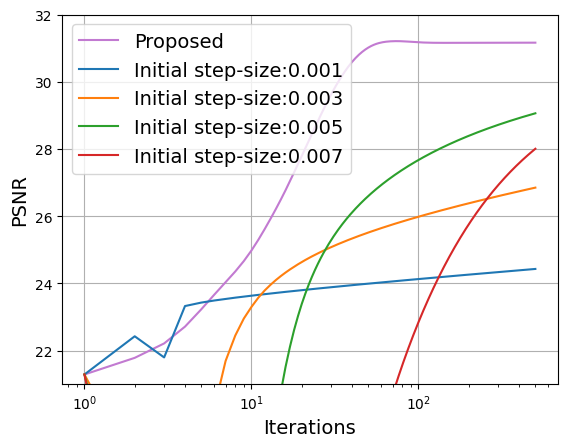}
     \centerline{SM-D}
   \end{minipage}    
\caption{CT: Comparison to subgradient methods. Note that SM-C is faster than SM-D in this case, and at times even faster than the proposed method initially.}
   \label{fig:ct_comp}
\end{figure}

\begin{figure}[!h]
\begin{minipage}{0.45\textwidth}
     \centering
     \includegraphics[width=\linewidth]{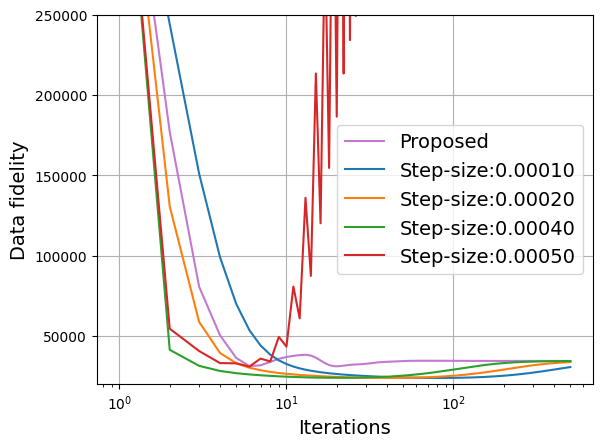}
   \end{minipage}
   \hspace{1cm}
    \begin{minipage}{0.45\textwidth}
     \centering
     \includegraphics[width=\linewidth]{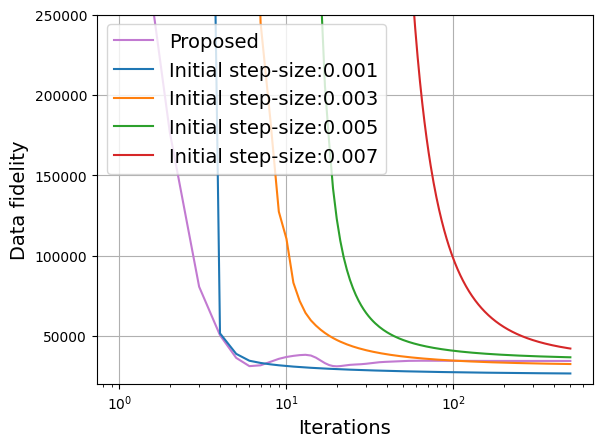}
   \end{minipage}    

\begin{minipage}{0.45\textwidth}
    \centering
     \includegraphics[width=\linewidth]{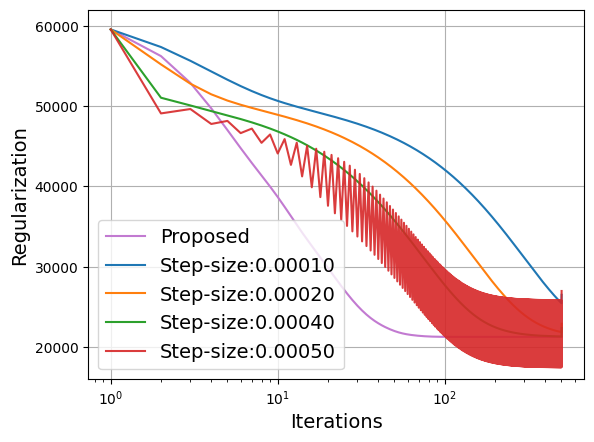}
     \hspace{1cm}
     \centerline{SM-C}
   \end{minipage}
   \hspace{1cm}
      \begin{minipage}{0.45\textwidth}
     \centering
     \includegraphics[width=\linewidth]{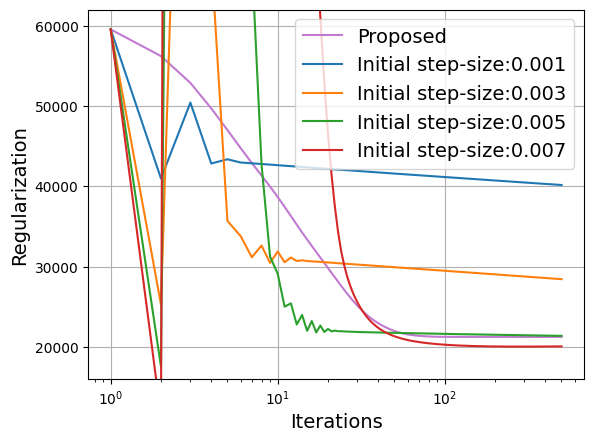}
     \centerline{SM-D}
   \end{minipage}    
\caption{CT: Data fidelity and regularization versus iterations plots. Notably, the subgradient methods with large step sizes exhibit oscillatory behavior, while the proposed method demonstrates more stable convergence.}
   \label{fig:ct_comp_dr}
\end{figure}

\textbf{Parameters:} We set $\gamma=600$. The parameters $c_0,c_1,c_2$ for the proposed method are chosen to be $500,100,0.1$. For SM-C, we select step-sizes from $\{1e$-$3$,$2e$-$3$,$4e$-$3$,$5e$-$3\}$, and $\{1e$-$2$,$3e$-$2$,$5e$-$2$,$7e$-$2\}$ as initial step-sizes for SM-D.

\textbf{Results:} Figure \ref{fig:ct_comp} compares the energy and PSNR plots of the proposed method and subgradient methods. While the constant step-size subgradient methods show substantial progress in the early iterations, they are quickly surpassed by the proposed method, which demonstrates a more consistent convergence. In contrast, the diminishing step-size subgradient methods exhibit much slower convergence overall. 

\begin{figure}[!h]
   \begin{minipage}{0.24\textwidth}
     \centering         
    \begin{tikzpicture}
    \node[anchor=south west,inner sep=0] (image) at (0,0) {\includegraphics[width=\textwidth]{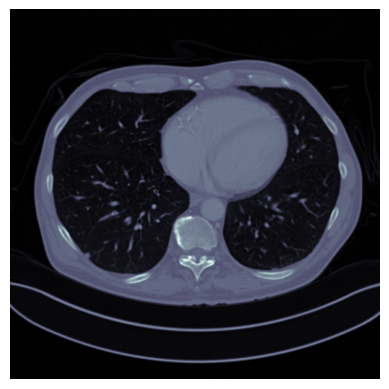}};
    \begin{scope}[x={(image.south east)},y={(image.north west)}]
    \node[text=white] at (0.85,0.9) {\footnotesize };
    \end{scope}
\end{tikzpicture}
     \centerline{Ground Truth}
   \end{minipage}
   \begin{minipage}{0.24\textwidth}
    \centering
    \begin{tikzpicture}
    \node[anchor=south west,inner sep=0] (image) at (0,0) {\includegraphics[width=\textwidth]{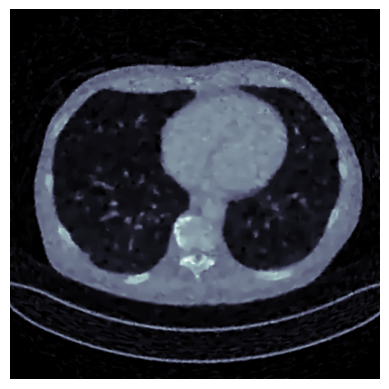}};
    \begin{scope}[x={(image.south east)},y={(image.north west)}]
    \node[text=white] at (0.85,0.9) {\footnotesize \textbf{31.06}};
    \end{scope}
\end{tikzpicture}
     \centerline{Proposed, 50 iter}
   \end{minipage}
   \begin{minipage}{0.24\textwidth}
    \centering
    \begin{tikzpicture}
    \node[anchor=south west,inner sep=0] (image) at (0,0) {\includegraphics[width=\textwidth]{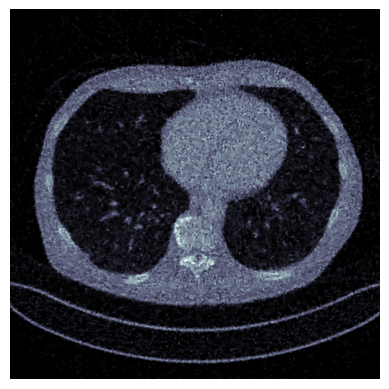}};
    \begin{scope}[x={(image.south east)},y={(image.north west)}]
    \node[text=white] at (0.85,0.9) {\footnotesize \textbf{26.19}};
    \end{scope}
\end{tikzpicture}
     \centerline{SM-C, 50 iter}
   \end{minipage}
   \begin{minipage}{0.24\textwidth}
    \centering
    \begin{tikzpicture}
    \node[anchor=south west,inner sep=0] (image) at (0,0) {\includegraphics[width=\textwidth]{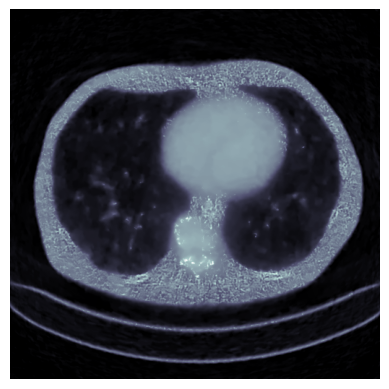}};
    \begin{scope}[x={(image.south east)},y={(image.north west)}]
    \node[text=white] at (0.85,0.9) {\footnotesize \textbf{26.67}};
    \end{scope}
\end{tikzpicture}
     \centerline{SM-D, 50 iter}
   \end{minipage}
   
    \begin{minipage}{0.24\textwidth}
     \centering
    \begin{tikzpicture}
    \node[anchor=south west,inner sep=0] (image) at (0,0) {\includegraphics[width=\textwidth]{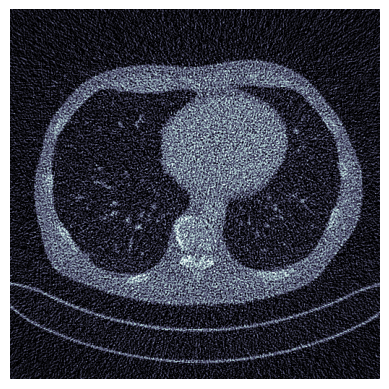}};
    \begin{scope}[x={(image.south east)},y={(image.north west)}]
    \node[text=white] at (0.85,0.9) {\footnotesize \textbf{21.29}};
    \end{scope}
\end{tikzpicture}
     \centerline{FBP}
   \end{minipage}
      \begin{minipage}{0.24\textwidth}
     \centering
    \begin{tikzpicture}
    \node[anchor=south west,inner sep=0] (image) at (0,0) {\includegraphics[width=\textwidth]{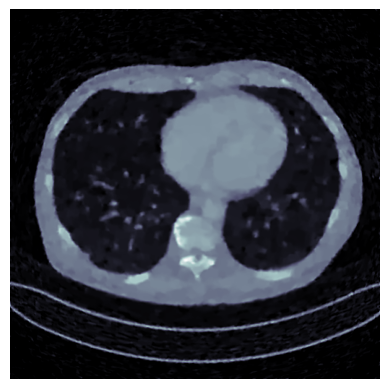}};
    \begin{scope}[x={(image.south east)},y={(image.north west)}]
    \node[text=white] at (0.85,0.9) {\footnotesize \textbf{31.16}};
    \end{scope}
\end{tikzpicture}
     \centerline{Proposed, 500 iter}
   \end{minipage}
      \begin{minipage}{0.24\textwidth}
     \centering
    \begin{tikzpicture}
    \node[anchor=south west,inner sep=0] (image) at (0,0) {\includegraphics[width=\textwidth]{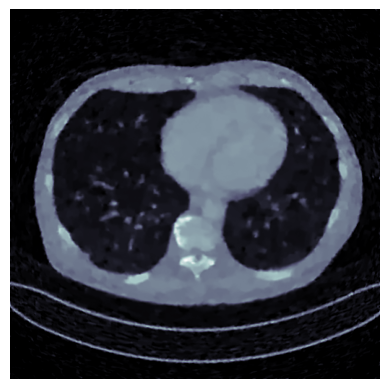}};
    \begin{scope}[x={(image.south east)},y={(image.north west)}]
    \node[text=white] at (0.85,0.9) {\footnotesize \textbf{31.16}};
    \end{scope}
\end{tikzpicture}
     \centerline{SM-C, 500 iter}
   \end{minipage}
      \begin{minipage}{0.24\textwidth}
     \centering
    \begin{tikzpicture}
    \node[anchor=south west,inner sep=0] (image) at (0,0) {\includegraphics[width=\textwidth]{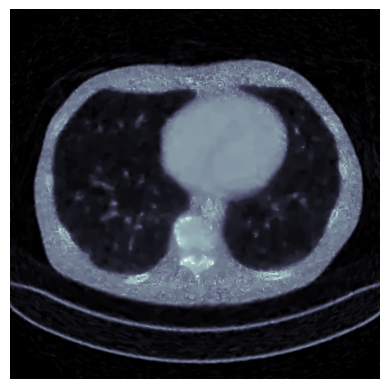}};
    \begin{scope}[x={(image.south east)},y={(image.north west)}]
    \node[text=white] at (0.85,0.9) {\footnotesize \textbf{29.06}};
    \end{scope}
\end{tikzpicture}
     \centerline{SM-D, 500 iter}
   \end{minipage}
   \caption{CT: Visual comparison of reconstructions, with PSNR shown at top right corner.}
   \label{fig:ct_recon}
\end{figure}

Additionally, Figure \ref{fig:ct_comp_dr} shows comparisons of the data fidelity and regularization term plots. All methods handle the data fidelity term reasonably well, though constant step-size subgradient methods sometimes drive it to values far below the eventual optimum, which may explain their faster initial decrease. In contrast, their behavior on the non-smooth regularization term differs. The constant step-size variants reduce the regularization term much more slowly than the proposed method, while the diminishing step-size variants exhibit oscillations in the early iterations. This highlights the superior stability of the proposed method throughout the optimization process. \rv{Interestingly, SM-C reaches a lower objective value than the proposed method after about 10 iterations, yet its PSNR is lower. A closer look at the plots reveals that SM-C overshoots the data fidelity term while progressing slowly on the regularizer, a consequence of using a single step-size, hence coupling both terms. In contrast, the proposed method balances the two terms more effectively.} Figure \ref{fig:ct_recon} shows the reconstructions at 50 and 500 iterations, further illustrating the effectiveness of the proposed method in producing high-quality results consistently.

\subsubsection{Deeper network}
\rv{We now consider a deeper ICNN with two more intermediate layers, which is represented by $R_{\boldsymbol{\theta}}(\mathbf{x})=\mathbf{W}_4h_4(\mathbf{W}_3\mathbf{P}\mathbf{z}_3+\mathbf{b}_3)$ with $\mathbf{z}_{i+1}=h_{i+1}(\mathbf{W}_i\mathbf{z}_i+\mathbf{b}_i)$ for $i=1,2$, and $\mathbf{z}_1=h_1(\mathbf{V}_0x+\mathbf{b}_0)$. Here $\mathbf{V}_0$, and $\mathbf{W}_i$ for $i=1,2$ correspond to convolution operators with 32 $5\times5$ output channels with 1 and 32 input channels respectively, and $\mathbf{P}$ denotes an average pooling operater with $16\times 16$ pool size. The fully connected layers $\mathbf{W}_3,\mathbf{W}_4$ consists of 256 and 1 output neurons respectively. The activations $h_i$ for $i=1,2,3$, and $h_4$ are chosen to be leaky ReLU and ReLU respectively, with the leaky ReLU's negative slope set to $0.2$. The regularizer is trained using the same data for 20 epochs with a learning rate of of $2e$-$5$, $(\beta_1,\beta_2)=(0.5,0.99)$, batch size of 20, and $\lambda_{GP}=20$. The corresponding reformulated problem can be written as:}
\rv{
\begin{equation}
    \min_{\mathbf{x,z}}D(\mathbf{Ax,y})+\gamma \mathbf{W}_4h_4(\mathbf{W}_3\mathbf{Pz}_3+\mathbf{b}_3)+\sum_{i=1}^2\delta_{C_{i+1}}(\mathbf{W}_i\mathbf{z}_i+\mathbf{b}_i,\mathbf{z}_{i+1})+\delta_{C_i}(\mathbf{V}_0\mathbf{x}+\mathbf{b}_0,\mathbf{z}_1)+\delta_{[0,\infty)}(\mathbf{x}).
\end{equation}
}
\rv{As the network depth increases, the reformulation only introduces additional indicator function for each intermediate layer, while the overall structure of the problem remains unchanged. The extra updates in the primal and dual variables correspond to analogous affine transformations and projection steps as in the shallower ones. Moreover, these operations can all be executed in parallel.} 

\rv{The choice of step-sizes can be naturally extended to deeper networks. Specifically, we pick step-sizes $\mathbf{T}=\operatorname{diag}(\tau_1\mathbf{I}_\mathbf{x},\tau_2\mathbf{I}_{\mathbf{z}_1},\tau_3\mathbf{I}_{\mathbf{z}_2},\tau_4\mathbf{I}_{\mathbf{z}_4}),\mathbf{S}=\operatorname{diag}(\sigma_0\mathbf{I}_{\mathbf{v}_0},\sigma_1\mathbf{I}_{\mathbf{v}_1},\sigma_2\mathbf{I}_{\mathbf{v}_2},\sigma_3\mathbf{I}_{\mathbf{v}_3},\sigma_4\mathbf{I}_{\mathbf{v}_4})$ given by:}
\rv{
\begin{equation}\label{step-size-ct-deep}
    \begin{aligned}
        \sigma_0=\frac{c_0}{\|\mathbf{A}\|^2},\sigma_1=\frac{c_1}{\|\mathbf{V}_0\|^2},\sigma_2=\frac{c_2}{\|\mathbf{W}_1\|^2},\sigma_3=\frac{c_3}{\|\mathbf{W}_2\|^2},\sigma_4=\frac{c_4}{\|\mathbf{W}_3\mathbf{P}\|^2},\\
    \tau_1=\frac{1}{c_0+c_1+c_2+c_3},\tau_2=\frac{1}{\sigma_1+c_2},\tau_3=\frac{1}{\sigma_2+c_3},\tau_4=\frac{1}{\sigma_3+c_4}.
    \end{aligned}
\end{equation} 
}

\begin{figure}[h!]
\begin{minipage}{0.45\textwidth}
     \centering
     \includegraphics[width=\linewidth]{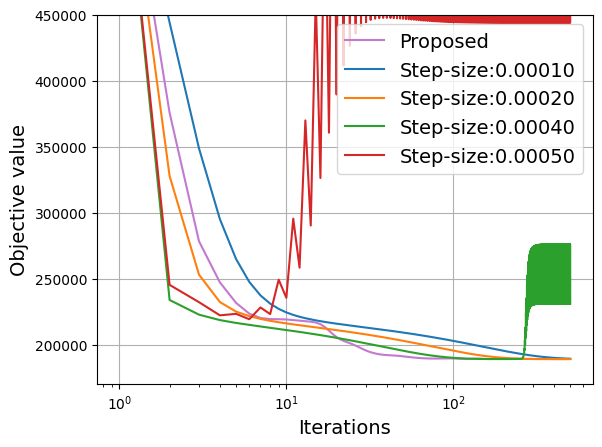}
   \end{minipage}
   \hspace{1cm}
    \begin{minipage}{0.45\textwidth}
     \centering
     \includegraphics[width=\linewidth]{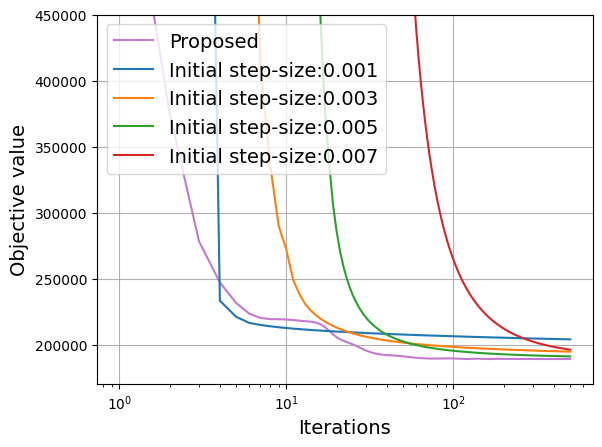}
   \end{minipage}    

\begin{minipage}{0.45\textwidth}
    \centering
     \includegraphics[width=\linewidth]{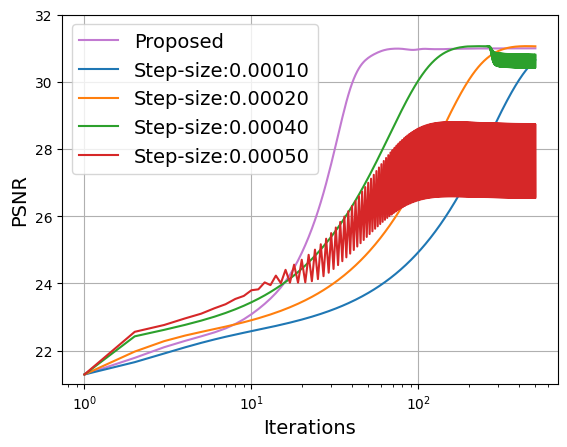}
     \hspace{1cm}
     \centerline{SM-C}
   \end{minipage}
   \hspace{1cm}
      \begin{minipage}{0.45\textwidth}
     \centering
     \includegraphics[width=\linewidth]{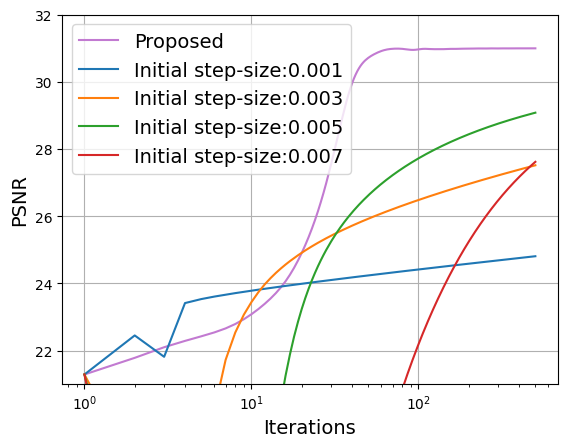}
     \centerline{SM-D}
   \end{minipage}    
\caption{CT: Objective and PSNR plots for deeper ICNN.}
   \label{fig:ct_deep_comp}
\end{figure}

\rv{
\textbf{Parameters:} We set $\gamma=650$. The parameters $c_0,c_1,c_2,c_3,c_4$ for the proposed method are chosen to be $500,100,5,0.5,0.1$. For SM-C, we select step-sizes from $\{1e$-$3$,$2e$-$3$,$4e$-$3$,$5e$-$3\}$, and $\{1e$-$2$,$3e$-$2$,$5e$-$2$,$7e$-$2\}$ as initial step-sizes for SM-D.
}

\rv{\textbf{Results:} Figure \ref{fig:ct_deep_comp} presents the comparison of energy and PSNR curves for the proposed method and subgradient methods when applied with a deeper ICNN. The overall behavior of the methods remains consistent with the shallower case. However, an important observation is that increasing the network depth does not lead to improvements in reconstruction quality in terms of PSNR. In fact, the deeper model performs on par with, or in some cases slightly worse than the shallower one. This suggests that simply increasing the network depth is not sufficient in this setting, and that the current adversarial training strategy may not fully exploit the potential of deeper architectures. This points to the need for exploring alternative training methods.}

\subsection{Image inpainting}
\label{sec:inpaint}
We consider an image inpainting task in this section. We randomly remove $30\%$ of the pixels of the image. We further add $3\%$ of Gaussian noise to the masked image. \rv{To demonstrate that our method performs well even for smooth regularizers, we consider a smooth version here. This also allows us to compare with more advanced optimization methods, such as accelerated algorithms. We utilize the same baseline architecture as in previous sections, while replacing the non-smooth activations with smooth functions. We consider the following smoothed approximation of ReLU:}
\rv{
\begin{equation*}
    \tilde{\psi}_{\nu}(\mathbf{x})=\begin{cases}
        0 &\text{ if }\mathbf{x}\leq0,\\
        \frac{\mathbf{x}^2}{2\nu}, &\text{ if }0<\mathbf{x}<\nu,\\
        \mathbf{x}-\frac{\nu}{2}, &\text{ otherwise},
    \end{cases}
\end{equation*}
}
\rv{where $\nu$ denotes a smoothing parameter. 
Similarly, we consider the smoothed approximation to leaky ReLU given by $\tilde{h}_1(\mathbf{x})=\kappa \mathbf{x}+(1-\kappa)\tilde{\psi}_{\nu}(\mathbf{x})$, where $\kappa$ corresponds to the negative slope of leaky ReLU.
Given the noise model, we adopt a $L^2$ data term and formulate the optimization problem as:}
\rv{
\begin{equation}
    \min_{\mathbf{x,z}}\frac{1}{2}\|\mathbf{Ax-y}\|^2_2+\gamma \mathbf{W}_2\tilde{h}_2(\mathbf{W}_1\mathbf{Pz}+\mathbf{b}_1)+\delta_{\tilde{C}_1}(\mathbf{V}_0\mathbf{x}+\mathbf{b}_0,\mathbf{z}),
\end{equation}
}
\rv{where $\mathbf{A}$ is a binary diagonal matrix that corresponds to the sampling mask. Here $\tilde{h}_2=\tilde{\psi}_{\nu}$ and $\tilde{C}_1=\{(p,q)|\tilde{h}_1(p)\leq q\}$. We apply the same dataset as in the salt and pepper experiment with the masked noisy images as degraded samples. The regularizer is trained for 20 epochs with a learning rate of $5e$-$4$, $(\beta_1,\beta_2)=(0.5,0.99)$, batch size of 8, and $\lambda_{GP}=5$. The updates of the primal-dual framework are as in (\ref{pdhg-l1}), with the $L^1$ data term replaced by the $L^2$ data term. The step-sizes are also chosen following (\ref{step-size}). Given the smoothness of the regularizer, we also consider the accelerated gradient method NMAPG \cite[Alg 4, supplementary material]{li2015accelerated} with a backtracking linesearch.} 

\rv{\textbf{Parameters:} We set $\gamma=0.1$. The parameters $c_1,c_2$ are picked as $0.01, 0.001$. For the gradient methods, we picked constant step-size $1.5$ and initial step-size $50$ for SM-C and SM-D respectively. The initial step-size for NMAPG is chosen to be $0.1$.}

\begin{table}[h!]
\centering
\rv{
\caption{Mean reconstructions time per image in seconds.}
\label{tab:ip-time}
\begin{tabular}{ccc}
Methods & Time (Mean$\pm$Std) & Speedup \\ \hline
Proposed & 0.18$\pm$0.026 & --- \\ 
SM-C & 0.27$\pm$0.047 & 1.52 \\ 
SM-D & 1.35$\pm$0.260 & 7.50 \\
NMAPG & 0.26$\pm$0.032 & 1.47
\end{tabular}%
}
\end{table}

\begin{figure}[!h]
\begin{minipage}{0.45\textwidth}
     \centering
     \includegraphics[width=\linewidth]{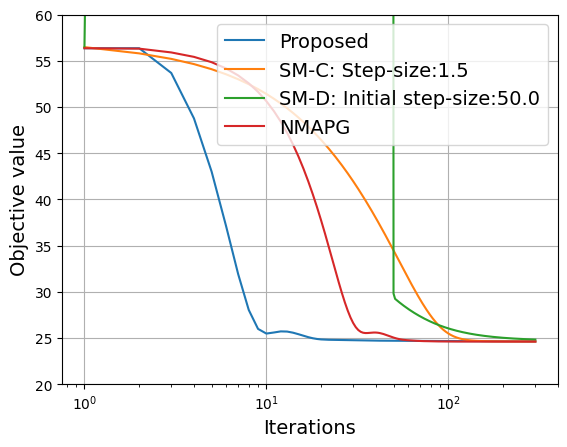}
   \end{minipage}
   \hspace{1cm}
    \begin{minipage}{0.45\textwidth}
     \centering
     \includegraphics[width=\linewidth]{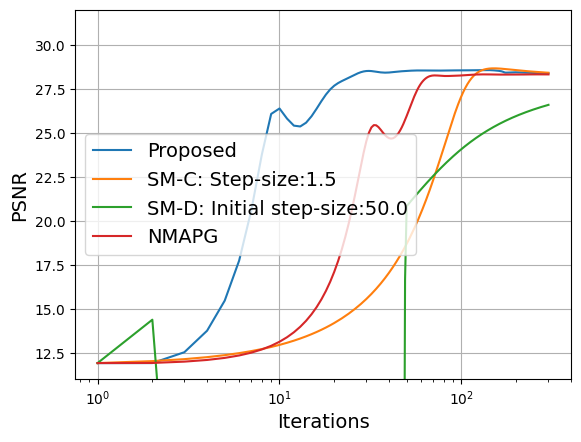}
   \end{minipage}    
\caption{Inpainting: Comparison to gradient-based methods.}
   \label{fig:ip_comp}
\end{figure}

\rv{
\textbf{Results:} To evaluate the performance of the proposed method across different test images, we solved the minimization problem on 20 test images and recorded the \rv{time} required to reduce the relative objective error below $1e$-$3$. Table \ref{tab:ip-time} shows the mean and standard deviation of the \rv{time in seconds}. Additionally, the mean speedup of the proposed method compared to the compared methods is reported, demonstrating the efficiency of the proposed method. Figure \ref{fig:ip_comp} show the comparisons of energy and PSNR plots. Notably, the proposed method converges significantly faster compared to both subgradient approaches. While the accelerated method outperforms the subgradient methods, it is still slower than the proposed approach.}

\subsubsection{Bilevel Training}
\rv{In this section, we investigate an alternative strategy for training the regularizer by incorporating the proposed method into the training pipeline via a bilevel optimization framework. The lower-level problem corresponds to the variational problem. We further absorb the regularization parameter inside the regularizer. For the upper-level loss, we consider the $\ell_2$ loss. The overall problem is formulated as:} 
\rv{\begin{equation}\label{eq:bilv-ip}
    \begin{aligned}
        &\min_{\boldsymbol{\theta}} \left\{\mathcal{L}(\boldsymbol{\theta}) = \frac{1}{N}\sum_{i=1}^N
    \|\hat{\mathbf{x}}_{\mathbf{y}_i}(\boldsymbol{\theta})-\mathbf{y}_i\|^2_2\right\}\\
    &\hat{\mathbf{x}}_{\mathbf{y}_i}(\boldsymbol{\theta}) = \arg\min\limits_{\mathbf{x}}\left\{\mathcal{J}_{\mathbf{y}_i}(\mathbf{x}; \boldsymbol{\theta}) = \frac{1}{2}\|\mathbf{Ax-y}_i\|^2_2 + R_{\boldsymbol{\theta}}(\mathbf{x}) \right\},
    \end{aligned}
\end{equation}
}
\rv{where $R_{\boldsymbol{\theta}}$ has the same architecture as the baseline one.}
\rv{In our experiments, we adopt the JFB approach for hypergradient computation. Specifically, after computing an approximate lower-level solution, we perform one additional gradient descent step and backpropagate through it. This is completely independent of the lower-level solver. We implement the bilevel learning using both the proposed solver and NMAPG as lower-level optimizers.}

\begin{table}[h!]
\centering
\rv{
\caption{Mean PSNR of reconstruction with different solvers and training schemes. The training time in seconds is shown in parentheses.}
\label{tab:ip-train-comp}
\begin{tabular}{cl cccc}
\toprule
& &  \multicolumn{4}{c}{Numerical Solver} \\
& & Proposed & SM-C & SM-D & NMAPG\\ 
\midrule 
\multirow{3}{*}{\rotatebox{90}{Training}\rotatebox{90}{Scheme}} 
& Adversarial (173) & 28.18 & 28.22 & 25.88 & 28.14\\ 
\cmidrule{2-6}
& Bilevel+NMAPG (435) & 29.46 & 29.40 & 27.31 
& 29.46 \\
& Bilevel+Proposed (431) 
& 29.80
& 29.72 
& 26.47 
& 29.80 \\ 
\bottomrule
\end{tabular}
}
\end{table}

\begin{figure}[!h]
   \begin{minipage}{0.19\textwidth}
     \centering         
    \begin{tikzpicture}
    \node[anchor=south west,inner sep=0] (image) at (0,0) {\includegraphics[width=\textwidth]{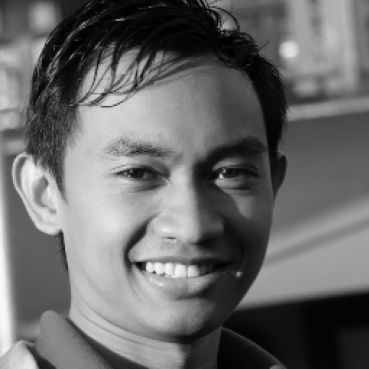}};
    \begin{scope}[x={(image.south east)},y={(image.north west)}]
    \node[text=white] at (0.80,0.9) {\footnotesize };
    \end{scope}
\end{tikzpicture}
    \centerline{Ground Truth}
   \end{minipage}
   \begin{minipage}{0.19\textwidth}
    \centering
    \begin{tikzpicture}
    \node[anchor=south west,inner sep=0] (image) at (0,0) {\includegraphics[width=\textwidth]{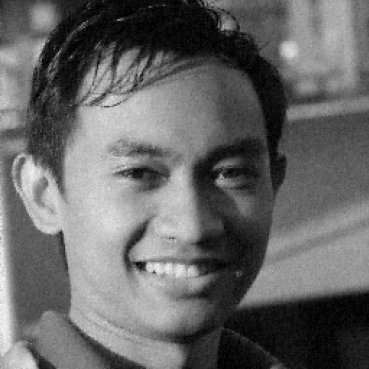}};
    \begin{scope}[x={(image.south east)},y={(image.north west)}]
    \node[text=white] at (0.80,0.9) {\footnotesize \textbf{29.93}};
    \end{scope}
\end{tikzpicture}
     \centerline{Proposed, 60 iter}
   \end{minipage}
   \begin{minipage}{0.19\textwidth}
    \centering
    \begin{tikzpicture}
    \node[anchor=south west,inner sep=0] (image) at (0,0) {\includegraphics[width=\textwidth]{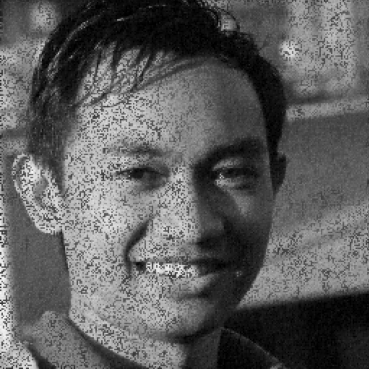}};
    \begin{scope}[x={(image.south east)},y={(image.north west)}]
    \node[text=white] at (0.80,0.9) {\footnotesize \textbf{19.54}};
    \end{scope}
\end{tikzpicture}
     \centerline{SM-C, 60 iter}
   \end{minipage}
   \begin{minipage}{0.19\textwidth}
    \centering
    \begin{tikzpicture}
    \node[anchor=south west,inner sep=0] (image) at (0,0) {\includegraphics[width=\textwidth]{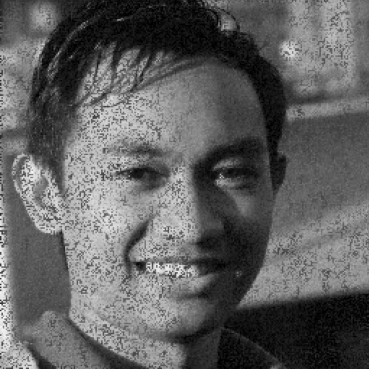}};
    \begin{scope}[x={(image.south east)},y={(image.north west)}]
    \node[text=white] at (0.80,0.9) {\footnotesize \textbf{20.82}};
    \end{scope}
\end{tikzpicture}
     \centerline{SM-D, 60 iter}
   \end{minipage}
   \begin{minipage}{0.19\textwidth}
     \centering
    \begin{tikzpicture}
    \node[anchor=south west,inner sep=0] (image) at (0,0) {\includegraphics[width=\textwidth]{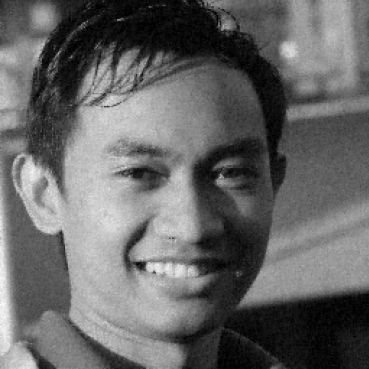}};
    \begin{scope}[x={(image.south east)},y={(image.north west)}]
    \node[text=white] at (0.80,0.9) {\footnotesize \textbf{29.75}};
    \end{scope}
\end{tikzpicture}
     \centerline{NMAPG, 60 iter}
   \end{minipage}%
   
    \begin{minipage}{0.19\textwidth}
     \centering
    \begin{tikzpicture}
    \node[anchor=south west,inner sep=0] (image) at (0,0) {\includegraphics[width=\textwidth]{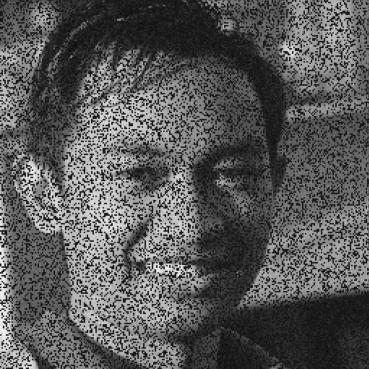}};
    \begin{scope}[x={(image.south east)},y={(image.north west)}]
    \node[text=white] at (0.80,0.9) {\footnotesize \textbf{11.91}};
    \end{scope}
\end{tikzpicture}
     \centerline{Noisy}
   \end{minipage}
    \begin{minipage}{0.19\textwidth}
     \centering
    \begin{tikzpicture}
    \node[anchor=south west,inner sep=0] (image) at (0,0) {\includegraphics[width=\textwidth]{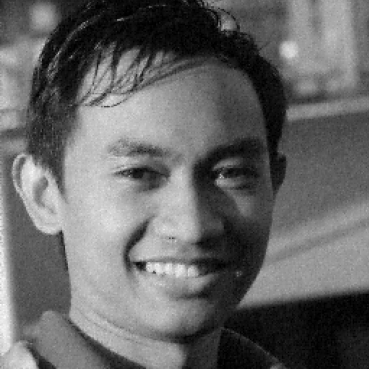}};
    \begin{scope}[x={(image.south east)},y={(image.north west)}]
    \node[text=white] at (0.80,0.9) {\footnotesize \textbf{30.01}};
    \end{scope}
\end{tikzpicture}
     \centerline{Proposed, 300 iter}
   \end{minipage}
    \begin{minipage}{0.19\textwidth}
     \centering
    \begin{tikzpicture}
    \node[anchor=south west,inner sep=0] (image) at (0,0) {\includegraphics[width=\textwidth]{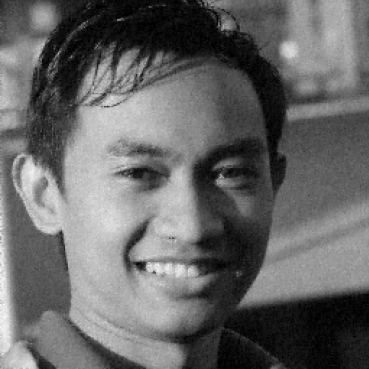}};
    \begin{scope}[x={(image.south east)},y={(image.north west)}]
    \node[text=white] at (0.80,0.9) {\footnotesize \textbf{29.97}};
    \end{scope}
\end{tikzpicture}
     \centerline{SM-C, 300 iter}
   \end{minipage}
    \begin{minipage}{0.19\textwidth}
     \centering
    \begin{tikzpicture}
    \node[anchor=south west,inner sep=0] (image) at (0,0) {\includegraphics[width=\textwidth]{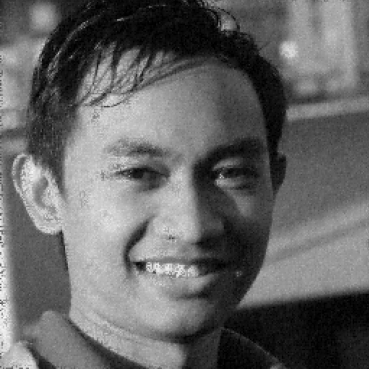}};
    \begin{scope}[x={(image.south east)},y={(image.north west)}]
    \node[text=white] at (0.80,0.9) {\footnotesize \textbf{27.85}};
    \end{scope}
\end{tikzpicture}
     \centerline{SM-D, 300 iter}
   \end{minipage}
   \begin{minipage}{0.19\textwidth}
     \centering
    \begin{tikzpicture}
    \node[anchor=south west,inner sep=0] (image) at (0,0) {\includegraphics[width=\textwidth]{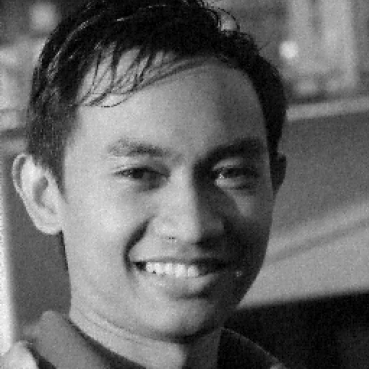}};
    \begin{scope}[x={(image.south east)},y={(image.north west)}]
    \node[text=white] at (0.80,0.9) {\footnotesize \textbf{30.01}};
    \end{scope}
\end{tikzpicture}
     \centerline{NMAPG, 300 iter}
   \end{minipage}
   \caption{Inpainting: Visual comparison of reconstructions, with PSNR shown at top right corner.}
   \label{fig:ip_recon_smooth}
\end{figure}

\rv{\textbf{Parameters:} We use 160 images and their corresponding noisy masked versions from the same training set as in the previous section. The norm of the gradient of the lower-level objective $\mathcal{J}$ serves as the stopping criterion, with tolerance set to $1e$-$1$. For the upper-level problem, we employ Adam with a learning rate $1e$-$3$ and $(\beta_1,\beta_2)=(0.9,0.99)$, and batch size of 4. Training is run for 30 epochs.}

\rv{\textbf{Results:} The learned regularizers are evaluated on the same 20 test images used in the previous experiments. We report the mean PSNR and training time for each lower-level solver, along with the corresponding results for the adversarially trained regularizer. Bilevel learning consistently improves reconstruction quality compared to adversarial training. Reconstructions obtained using the regularizer learned in the bilevel setting, with the proposed method as the lower-level solver, are shown in Figure \ref{fig:ip_recon_smooth}. The proposed method outperforms standard gradient-based methods and shows a slight advantage over NMAPG during the early iterations. Figure \ref{fig:ip_comp_train} compares reconstructions from regularizers trained with different methods, demonstrating that bilevel learning produces sharper results.}

{
\newcommand{\addfigure}[3]{%
    \begin{tikzpicture}%
    \node[anchor=south west,inner sep=0] (image) at (0,0) {\includegraphics[width=.6\textwidth]{images/inpaint/#2.png}};%
    \begin{scope}[x={(image.south east)},y={(image.north west)}]%
    \node[text=white, anchor=north east] at (0.99,0.99) {\footnotesize \textbf{#3}};
    \node[text=white, anchor=south west] at (0.01,0.01) {\footnotesize \textbf{#1}};%
    \end{scope}%
    \end{tikzpicture}%
}
\newcommand{\addcrop}[3]{%
    \begin{tikzpicture}%
    \node[anchor=south west,inner sep=0] (image) at (0,0) {\includegraphics[width=\textwidth,trim={90 80 90 110}, clip]{images/inpaint/#2.png}};%
    \begin{scope}[x={(image.south east)},y={(image.north west)}]%
    \node[text=white, anchor=north east] at (0.99,0.99) {\footnotesize \textbf{#3}};
    \node[text=white, anchor=south west] at (0.01,0.01) {\footnotesize \textbf{#1}};%
    \end{scope}%
    \end{tikzpicture}%
}
\begin{figure}[!h]
    
   \begin{minipage}{0.19\textwidth}
       \addcrop{Ground Truth}{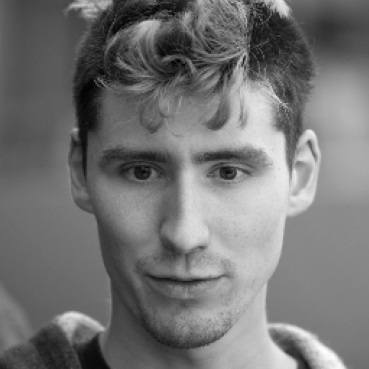}{ }
   \end{minipage}
    \begin{minipage}{0.19\textwidth}
    \addcrop{Noisy}{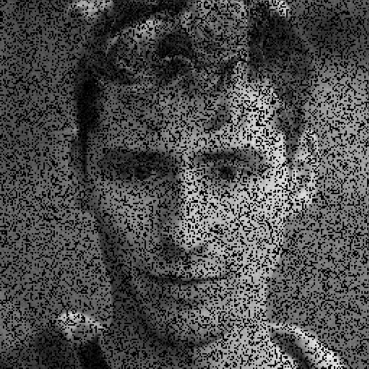}{12.21}
    \end{minipage}
    \begin{minipage}{0.19\textwidth}
        \addcrop{AR}{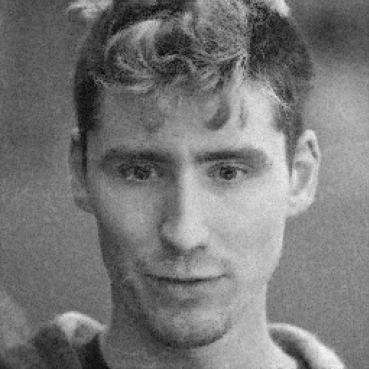}{28.38}
    \end{minipage}
    \begin{minipage}{0.19\textwidth}
    \addcrop{BL+NMAPG}{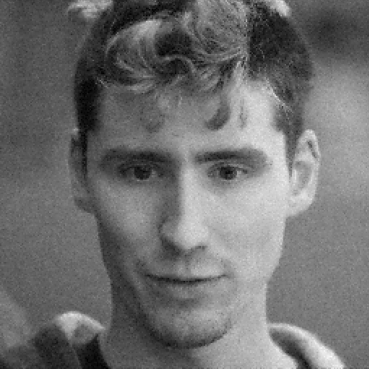}{30.27}
    \end{minipage}
    \begin{minipage}{0.19\textwidth}
    \addcrop{BL+Proposed}{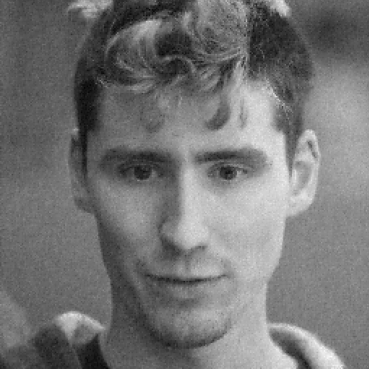}{30.47}
    \end{minipage}
    \caption{Inpainting: Visual comparison of reconstructions with different training scheme. }
   \label{fig:ip_comp_train}
\end{figure}
}

\section{Conclusion}
\label{sec:conclusion}
We proposed an efficient method for solving the optimization problem in variational reconstruction with a learned convex regularizer. A key challenge comes from the non-smoothness of the ICNN regularizer, whose proximal operator lacks a closed-form solution. To overcome this, we decoupled the neural network layers by introducing auxiliary variables corresponding to the layer-wise activations. While this initially resulted in a non-convex problem, we drew inspiration from the convexity of epigraphs and reformulated it as a convex optimization problem. We then proved that this reformulation is equivalent to the original variational problem and applied a primal-dual algorithm to solve it. 

Numerical experiments demonstrated that the proposed method not only outperforms subgradient methods in terms of convergence speed but also exhibits greater stability throughout the optimization process, as evidenced by the smoother behavior observed in the energy versus iterations plots, as well as those depicting data fidelity and regularization. \rv{Moreover, in the smooth setting, our approach compares favorably even against accelerated methods, highlighting its effectiveness not just in the non-smooth setting. We further showed that the method can be integrated into a bilevel training framework, where it serves as a lower-level solver for learning the regularizer. This integration yielded improvements in reconstruction quality compared to adversarially trained regularizers, demonstrating the significance of the proposed method not only for reconstruction but also for training.}

Additionally, we note that the updates of the proposed method are independent, enabling parallel computation. Looking forward, we aim to explore the potential of extending the proposed method to primal-dual variants that leverage this, such as coordinate-descent primal-dual algorithms \cite{fercoq2019coordinate}.


\bibliographystyle{siamplain}
\bibliography{references}
\end{document}


\maketitle

\section{A detailed example}

Here we include some equations and theorem-like environments to show
how these are labeled in a supplement and can be referenced from the
main text.
Consider the following equation:
\begin{equation}
  \label{eq:suppa}
  a^2 + b^2 = c^2.
\end{equation}
You can also reference equations such as \cref{eq:matrices,eq:bb} 
from the main article in this supplement.

\lipsum[100-101]

\begin{theorem}
An example theorem.
\end{theorem}

\lipsum[102]
 
\begin{lemma}
An example lemma.
\end{lemma}

\lipsum[103-105]

Here is an example citation: \cite{KoMa14}.

\section[Proof of Thm]{Proof of \cref{thm:bigthm}}
\label{sec:proof}

\lipsum[106-112]

\section{Additional experimental results}
\Cref{tab:foo} shows additional
supporting evidence. 

\begin{table}[htbp]
\footnotesize
  \caption{Example table.}  \label{tab:smfoo}
\begin{center}
  \begin{tabular}{|c|c|c|} \hline
   Species & \bf Mean & \bf Std.~Dev. \\ \hline
    1 & 3.4 & 1.2 \\
    2 & 5.4 & 0.6 \\ \hline
  \end{tabular}
\end{center}
\end{table}

\bibliographystyle{siamplain}
\bibliography{references}